\documentclass[11pt, reqno]{amsart}
\usepackage{amssymb,amsmath}

 \usepackage{fancyhdr} 

\theoremstyle{plain}
 \newtheorem{thm}{Theorem}
 \newtheorem{cor}[thm]{Corollary}
 \newtheorem{lem}[thm]{Lemma}
 \newtheorem{prp}[thm]{Proposition}

\def\thetitle{Inflectional loci of scrolls  II}
\title{Inflectional loci of scrolls II}
\author[A. Lanteri, R. Mallavibarrena, R. Piene]{Antonio Lanteri,
Raquel Mallavibarrena, Ragni Piene}
\address{Dipartimento di Matematica ``F. Enriques'',
Universit\`a degli Studi di Milano, Via C. Saldini, 50,  I-20133
Milano, Italy} \email{antonio.lanteri@unimi.it}
\address{Departamento de Algebra, Facultad de Ciencias Matem\`aticas,
Plaza de las Ciencias, 3 - Universidad Complutense de Madrid,
28040 Madrid, Spain}
\email{rmallavi@mat.ucm.es}
\address{CMA/Department of Mathematics, University of Oslo, P. O. Box 1053
Blindern, NO-0316 Oslo, Norway}
\email{ragnip@math.uio.no}

\date{\today}

\begin{document}

\begin{abstract}
Let $X\subset \mathbb P^N$ be a scroll over a $m$-dimensional
variety $Y$. We find the locally free sheaves on $X$ governing
the osculating behavior of $X$, and, under certain
dimension assumptions, we compute the cohomology class and the
degree of the inflectional locus of $X$. The case $m=1$ was
treated in \cite{LMP}. Here we treat the case $m\ge 2$, which is
more complicated for at least two reasons: the expression for the
osculating sheaves and the computations of the class of the
inflectional locus become more complex, and the dimension
requirements needed to ensure validity of the formulas are more
severe.
\end{abstract}

\maketitle

\section{Introduction}

Let $X \subset \mathbb{P}^N=\mathbb{P}(V)$ be a non-degenerate,
smooth, complex projective variety of dimension $n$, let
$\mathcal{L}$ denote the hyperplane bundle, and identify $V$ with
a subspace of $H^0(X,\mathcal{L})$. Let $\mathcal P_{X}^k(\mathcal
L)$ be the $k$th principal parts bundle of $\mathcal L$, and let
$j_k:V_X= V \otimes \mathcal{O}_X \to
\mathcal{P}^k_X(\mathcal{L})$ be the sheaf homomorphism
associating to every section $\sigma \in V$ its $k$th jet
evaluated at $x$, for every $x \in X$. The homomorphism $j_{k,x}$
allows us to define the $k$th osculating space to $X$ at $x$ as
follows:  $Osc^k_x(X):=\mathbb P({\rm Im}\,(j_{k,x}))$ (we
interpret $\mathbb P(V)$ as the
 set of codimension 1 vector subspaces of $V$). Let $s_k$ denote the
generic rank of $j_{k,x}$. The $k$th inflectional locus
$\Phi_k(X)$ of $X$ is defined to be the locus where the rank
of $j_{k,x}$ is strictly less than $s_k$.
 In \cite{LMP} we studied the particular case where $X$ is a scroll over a smooth curve.
In this case the expected generic rank of $j_k$ is $kn+1$, for
every $k$ such that $kn \leq N$. For the largest such $k$, we
obtained an explicit formula for the integral cohomology class of
$\Phi_k(X)$. In particular, it allowed us to count the flexes when
there are finitely many, and to characterize uninflected scrolls in an appropriate range.

When we consider scrolls $X$ over a higher dimensional smooth projective variety $Y$,
the situation becomes more difficult for several reasons. First of
all, taking as starting point the general framework of \cite{LMP},
the generalization of the construction of the osculating
bundles when $m :=\rm {dim}Y\geq 2$ is not straightforward and requires the
definition of a new bundle $\mathcal{T}_k$ to get the appropriate
diagrams to determine the inflectional locus. Secondly, it becomes more and more difficult to obtain
explicit expressions for the cohomology class of $\Phi_k(X)$, using
Porteous' formula, as the dimension
of $Y$ increases. Like in the case of scrolls over curves
\cite{LMP}, we need to impose that $k$ is such that $s_k$ is the
maximum value allowed for a scroll of dimension $n$ over $Y$ in
$\mathbb P^N$. This requirement on $m$, $n$, and $N$ fixes the
range of validity of the results so that there are many examples
that do not fulfill the requirement and that we therefore cannot treat
with the present methods. Some of these examples, also
interesting from the point of view of uninflectedness, are
discussed in this paper. The dimension requirements are explained
in Section 2, where we describe the general setting.

The main general results of the paper are Theorems \ref{thm-1} and
\ref{thm-2}. Theorem \ref{thm-1} provides the technical tools
leading to a relevant exact diagram. With Theorem \ref{thm-2} we
extract from this diagram a formal expression for the cohomology
class of $\Phi_k(X)$, which is given by the product of $k+1$
inverses of the total Chern classes of some vector bundles
involving symmetric powers of $T_Y$, the vector bundle
$\pi_*\mathcal{L}$, where $\pi:X \to Y$ is the scroll projection,
and $\mathcal{L}$ itself. In order to make this formula explicit
one needs either to simplify the vector bundles involved or to
reduce the number of factors. The first approach is developed in
Section 3, where we treat the case that $Y$ is an abelian variety.
Here the complexity of the computation is represented only by the
dimension $m$, and in this case we are able to get explicit
formulas for the class of $\Phi_k(X)$ when $m=2$ (Proposition
\ref{Theorem-abelian}) and $m=3$ (Proposition
\ref{Theorem-abelian3}). The second approach takes the rest of the
paper. Here the complexity depends on the codimension $\ell$ of
$\Phi_2(X)$. When $\ell=1$, in Section \ref{k=2,divisor} we give
an explicit formula for $\Phi_2(X)$ holding for any $m$ and $n$
(Theorem \ref{thm-divisor}), a lower bound for its degree and a
better lower bound holding if a suitable adjoint bundle on $Y$ is
nef. This situation is elaborated further in Section 5,
assuming moreover that $m=2$ and $n=3$. In this case, for
threefold scrolls in $\mathbb P^8$, Theorem \ref{thm-details}
gives the lower bound $\deg \Phi_2(X)\ge 3d$, where $d$ is the
degree of $X$, with a few explicit exceptions. In the last two
sections we take care of the case $\ell \geq 2$. In Section \ref{k=2}
we make explicit the expression of $\Phi_2(X)$ for $(n,m)=(3,2)$
and we discuss the uninflectedness of $X$ for $Y$ belonging to
various classes of surfaces. Finally, in Section \ref{k=2,m=3},
explicit formulas for $\Phi_2(X)$ are computed, for the various $\ell$, when $X$ is a four-dimensional scroll over a threefold $Y$, and
we prove that for $Y=\mathbb P^3$ or $\mathbb Q^3$ none of these
scrolls are uninflected.

\section{Scrolls}
Let $X \subset \mathbb P^N=\mathbb P(V)$ be a non-degenerate, smooth
$n$-dimensional scroll over a smooth $m$-dimensional variety $Y$, and let $\pi:X
\to Y$ be the projection. We identify $V$ with a vector subspace
of $H^0(X, \mathcal L)$, where $\mathcal L$ is the hyperplane line
bundle of $X$.  Since $X$ is smooth and $\mathcal L$ is very ample, the homomorphism
$j_{1,x}$ has rank $n+1$ at every point $x \in X$. We are
interested in the generic rank of the homomorphism $j_{k}$ for
$k \geq 2$. Choose local coordinates $(u_1,\dots, u_m, v_{m+1}, \dots ,v_n)$
around $x \in X$ in such a way that $(u_1,\dots, u_m)$ are (pullbacks of)  local
coordinates around $\pi(x)$ on $Y$, while $(v_{m+1}, \dots, v_n)$ are
local coordinates around $x$ on the fiber $\pi^{-1}(\pi(x))\cong \mathbb P^{n-m}$. Then, any section $\sigma \in V$ can be represented, locally
around $x$, in the following form
$$ \sigma= a(u_1,\dots, u_m) + \sum_{j=m+1}^n b_j(u_1,\dots, u_m) v_j.$$
Since $\sigma$ is linear in the variables $v_j$ for $j=m+1, \dots,
n$, when we take the derivatives of order $h$ for $2 \leq h \leq
k$, only the following are not obviously zero:
$\sigma_{u_1^{i_1},\cdots, u_m^{i_m}}$ with $\sum_{j=1}^m i_j=h$ and
$\sigma_{u_1^{i_1},\cdots, u_m^{i_m},v_j}$ with $\sum_{j=1}^m i_j=h-1$
and $j=m+1, \dots n$. Those of the former type
are counted by the forms of degree $h$ in $(u_1,\dots, u_m)$, hence their
number is $\binom{m-1+h}h$; those of the latter type, once $v_j$ is fixed,
are counted by the forms of degree $h-1$ in $(u_1,\dots, u_m)$, which are
$\binom{m-1+h-1}{h-1}$, hence their total number is $(n-m)\binom{m-1+h-1}{h-1}$. Let $s_k$ be the
generic rank of $j_{k}$. Then the above says that $s_k$ cannot
exceed
\[ \textstyle n+1 +\sum_{h=2}^k \bigl(\binom{m-1+h}h + (n-m)\binom{m-1+h-1}{h-1}\bigr)
\]

\begin{lem}
The generic rank $s_k$ of $j_{k}$ satisfies
\begin{equation} \label{tag 1}
s_k \leq r_k,
\end{equation}
where
\begin{equation} \label{tag 2}\textstyle
 r_k :=(n-m)\binom{m+k-1}{k-1} + \binom{m+k}k
\end{equation}
\end{lem}

\begin{proof}
Use the fact that
\[\textstyle \sum_{h=0}^k \binom{m-1+h}h=\binom{m+k}k.
\]
\end{proof}

\noindent{\it{Example 1.}} Let $Y \subset \mathbb{P}^M$ be a
$m$-dimensional smooth variety, let ${\mathcal H}$ be the
hyperplane line bundle, and consider the scroll
$X:=\mathbb{P}_Y({\mathcal H}^{\oplus (n-m+1)}) \subset
\mathbb{P}^N$, where $N=(M+1)(n-m+1)-1$. Clearly, $X \cong Y
\times \mathbb{P}^{n-m}$, embedded by the Segre embedding of
$\mathbb P^M\times \mathbb P^{n-m} \to \mathbb P^N$. Fix any point
$x \in X$ and let $y:=\pi(x)$, where $\pi:X \to Y$ is the
projection. Arguing as in \cite[p. 1049--1050]{PS} or \cite[Lemma
2.1]{LM}, one can express the matrix representing $j_{k,x}$ for $X
\subset \mathbb{P}^N$ in terms of the blocks representing
$j_{k,y}$ and $j_{k-1,y}$ for $Y \subset \mathbb{P}^M$, showing
that
$$\text{rk}\, j_{k,x} = (n-m) \text{rk}\, j_{k-1,y} +
\text{rk}\, j_{k,y}.$$ In particular, we see that in this special
case, equality holds in (\ref{tag 1}) if and only if the base
variety $Y \subset \mathbb{P}^M$ is generically $k$-regular.
Moreover, recalling that $\Phi_{k-1}(Y) \subseteq \Phi_k(Y)$, we
see that $\Phi_k(X) = \pi^{-1}(\Phi_k(Y))$. So, either $X$ is
uninflected, i.\,e., $\Phi_k(X) = \emptyset$, or $\dim \Phi_k(X) =
\dim \Phi_k(Y) + (n-m)$. In particular, $X$ is uninflected if and
only if $Y$ is.
\par

\medskip
\noindent{\it{Example 2.}} Here are some examples of scrolls for
which (\ref{tag 1}) is a strict inequality. First of all, we
mention the Segre products: a) $X= \mathbb{P}^2 \times
\mathbb{P}^1 \subset \mathbb{P}^5$ and b) $X= \mathbb{P}^1 \times
\mathbb{P}^1 \times \mathbb{P}^1 \subset \mathbb{P}^7$. According
to what we said in Example 1, they cannot satisfy the equality
$s_2=r_2$. This also follows from the fact that $r_2=9$, while the
dimension $N$ of the embedding space is too small in both cases.
Note that in the former case, $X$ is also a scroll over
$\mathbb{P}^1$, so that this situation is covered by \cite{LMP}.
Look at the latter case: since $X = \mathbb{P}_{\mathbb{P}^1
\times \mathbb{P}^1}(\mathcal{O}(1,1)^{\oplus 2})$, it follows
from Example 1 that $\text{rk}\, j_{2,x}=7$ for every $x \in X$. This
makes $X$ an example of a perfectly hypo-osculating threefold scroll,
but, unfortunately, we cannot include it in our subsequent
discussion.

More generally, we can consider the following Segre products.
\medskip

A) $X = \mathbb{P}^m \times \mathbb{P}^{n-m} \subset
\mathbb{P}^N$, where $N=m(n-m)+n$. By using coordinates $(1:u_1:
\dots :u_m)$ on $\mathbb{P}^m$ and $(1:v_1: \dots :v_{n-m})$ on
$\mathbb{P}^{n-m}$, $X$ is locally represented in $\mathbb{P}^N$
by coordinates
$$(1: \dots: u_i: \dots :v_j: \dots :u_iv_j: \dots ).$$
Clearly all second derivatives of such coordinates are zero and we
immediately see that $j_{2,x}$ has rank $1 + m + (n-m) + m(n-m)$
at every point $x \in X$. Thus $s_2=N+1$. On the other hand,
from (\ref{tag 2}) we get $r_2 = (n-m)(m+1)+
\binom{m+2}2$. Therefore $s_2 = r_2 - \binom{m+1}2$.

B) $X=(\mathbb{P}^1)^n \subset \mathbb{P}^{2^n-1}$. Clearly, $X =
\mathbb{P}_Y(\mathcal{O}_Y(1,1, \dots ,1)^{\oplus 2})$, where
$Y=(\mathbb{P}^1)^{n-1}$. Thus, if $u_i$ is a local parameter on
the $i$th factor of $Y$, $i=1, \dots, n-1$, $X$ is locally
represented in $\mathbb{P}^{2^n-1}$ by coordinates
$$ (1 : u_1: \dots:v:u_1u_2 : \dots :u_{n-1}v:\dots :u_1 \dots u_{n-1} v).$$ Since any such coordinate $\sigma$ is linear in the
parameters we have that in the matrix representing $j_{2,x}$ the
$n-1$ rows representing the second derivatives $\sigma_{u_j^2}$
for $j=1, \dots ,n-1$ (in addition to $\sigma_{v^2}$) are zero. The remaining rows are clearly linearly independent, hence
$s_2 =1+n+\binom{n}2= \binom{n+2}2 -n= r_2-(n-1)$. In fact, it is easy to see that \[\textstyle s_{n-1}=1+\binom{n}1 +\dots \binom{n}{n-1}=2^n-1<r_{n-1}=3\binom{2n-3}{n-1},\] and that $X$ is perfectly hypo-osculating.

Another interesting threefold scroll excluded from our discussion, that
we would like to mention, is $\mathbb{P}(T_{\mathbb{P}^2}) \subset
\mathbb{P}^7$, for which we have $s_2 =8 < 9=r_2$.

Next, let $Y (\cong \mathbb{F}_1) \subset \mathbb{P}^4$ be the
rational cubic surface scroll, and let $X = Y \times \mathbb{P}^1
\subset \mathbb{P}^9$, embedded via the Segre embedding. As $Y$ is
locally described by $(1:u_1:u_2:u_1u_2:u_1^2u_2)$, $X$ is locally
represented by
$$(1:u_1:u_2:v:u_1u_2:u_1v:u_2v:u_1u_2v:u_1^2u_2:u_1^2u_2v).$$
Note that no term of degree $\geq 2$ in $u_2$ appears. Hence $s_2
\leq 8$, and in fact a direct check shows that $s_2=8$. Thus $s_2
< r_2=9$.

\medskip

We want to investigate the inflectional loci of non-degenerate, smooth scrolls over a smooth base space.
We shall do this under the assumption
that equality holds in (\ref{tag 1}). Note that this is
a serious restriction, as the previous examples show.

\begin{thm}\label{thm-1} Let $X\subset \mathbb P^N$ be a $n$-dimensional scroll
over a smooth $m$-dimensional variety $Y$, with hyperplane bundle $\mathcal L=\mathcal O_{P^N}(1)|_X$.
For given $k\ge 1$, set  $r_k=(n-m)\binom{m+k-1}{m} + \binom{m+k}m $. Assume that
the generic rank of $j_k$ is equal to $r_k$ for all $k\ge 1$ such that $r_k \le N+1$.
Set $\mathcal Q_k={\rm Coker}\, j_k$. Then, for all such $k$,
\begin{itemize}
\item[(i)] $\mathcal Q_k^{\vee}$ is a locally free sheaf of rank
$\binom{n+k}n-r_k$; \item[(ii)] there exist locally free sheaves
$\mathcal M_k$, of rank $\binom{n+k-1}{n-1}-\binom{m+k-1}{m-1}-
(n-m)\binom{m+k-2}{m-1}$, and $\mathcal T_k$, of rank
$(n-m)\binom{m+k-2}{m-1}+\binom{m+k-1}{m-1}$,  and exact sequences
\begin{equation*}
0\to \mathcal Q_{k-1}^{\vee} \to \mathcal Q_k^{\vee} \to \mathcal M_k^\vee \to 0,
\end{equation*}
and
\begin{equation*}
0\to \mathcal T_k \to S^k\Omega_X\otimes L \to \mathcal M_k \to 0;
\end{equation*}
\item[(iii)] the quotient sheaves $\mathcal E_k:= \mathcal
P_X^k(\mathcal L)^{\vee}/\mathcal Q_k^{\vee}$ are locally free, of
rank $r_k$, and there exist exact sequences
\begin{equation*}
0\to \mathcal E_{k-1} \to \mathcal E_k \to \mathcal T_k^\vee \to 0.
\end{equation*}
\end{itemize}
\end{thm}

\begin{proof}
We will adapt the proof of \cite[Thm. 1]{LMP}.
Define the sheaf $\mathcal E_k$ by the exact sequence
\[0 \to \mathcal Q_k^{\vee} \to \mathcal P^k_X(\mathcal L)^{\vee}
\to \mathcal E_k \to
0.\]
Set $\mathcal K_k:\ = \text{ker}(S^{k-1}\Omega_Y
\otimes \Omega_Y \to S^k \Omega_Y)$; then we have an exact sequence
\begin{equation}\label{defK}
0 \to \mathcal K_k \to S^{k-1}\Omega_Y \otimes \Omega_Y \to
S^k\Omega_Y \to 0,
\end{equation}
and $\mathcal K_k$ is locally free, with rank
$m\binom{m+k-2}{m-1}-\binom{m+k-1}{m-1}$.
We observe that $\pi^*\mathcal K_k$ is the kernel of the composite
map
\begin{equation}\label{comp}
\pi^*S^{k-1}\Omega_Y \otimes \Omega_X \to S^{k-1}\Omega_X \otimes
\Omega_X \to S^k \Omega_X.
\end{equation}
To see this, let $x \in X$ and let $u_1,\ldots,u_m$ denote local
coordinates on the base variety $Y$ around $\pi(x)$, and $v_{m+1},
\dots ,v_n$ local coordinates on the fiber through $x$ around $x$.
So, if $A:= \mathcal O_{X,x}$, we have the following isomorphisms:
\setlength\arraycolsep{2pt}
\begin{eqnarray*} (\pi^*\Omega_Y)_x &\cong & \textstyle A\ du_1 \oplus \cdots \oplus A\ du_m,\\
 \Omega_{X,x} &\cong & \textstyle \bigoplus_{i=1}^m A\ du_i \oplus \bigoplus_{j=m+1}^n
A\ dv_j,\\
 (S^{k-1}\Omega_X)_x &\cong & \textstyle \bigoplus_{i_1, \dots , i_n}
A\ du_1^{i_1}\cdots du_m^{i_m} dv_{m+1}^{i_{m+1}} \cdots dv_n^{i_n},
\end{eqnarray*}
with
$\sum_{j} i_j = k-1$. The composite map (\ref{comp}) can be
described locally in this way:
$$du_1^{i_1}\cdots du_m^{i_m} \otimes du_j \mapsto du_1^{i_1}\cdots du_j^{i_j+1}\cdots du_m^{i_m},$$
for $j=1,\ldots,m$, and
$$du_1^{i_1}\cdots du_m^{i_m}  \otimes dv_h \mapsto du_1^{i_1}\cdots du_m^{i_m} dv_h$$
for $h=m+1,\ldots,n$.

The kernel is generated by the elements
$$du_1^{i_1}\cdots du_j^{i_j}\cdots du_l^{i_l-1}\cdots du_m^{i_m}\otimes du_l
-du_1^{i_1}\cdots du_j^{i_j-1}\cdots du_l^{i_l}\cdots du_m^{i_m}\otimes du_j, $$
where $ i_1+\cdots i_m=k$, $j\neq l$, and $i_j, i_l \ge 1$ These generators are not linearly independent when $m> 2$.
As $\pi^*\mathcal K_k$ is a
subbundle of $\pi^*S^{k-1}\Omega_Y \otimes \Omega_X$ via the map
$$0 \to \pi^*S^{k-1}\Omega_Y \otimes \pi^*\Omega_Y \to
\pi^*S^{k-1}\Omega_Y \otimes \Omega_X$$ it is obvious that it
coincides with the kernel just computed.

Moreover the composite map (\ref{comp}) has constant rank, so
its image is a subbundle of $S^k \Omega_X$. Local
computations show that the image in $S^k\Omega_X$ is generated by
$\{du_1^{i_1}\cdots du_m^{i_m}; du_1^{i'_1}\cdots du_m^{i'_m}dv_h\}$
with $\sum i_j=k$, $\sum i'_j=k-1$, $h=m+1, \dots ,n$.
These generators are independent and their
number is $\binom{m+k-1}{m-1}+(n-m)\binom{m+k-2}{m-1}$.
This is valid for all points $x
\in X$. Tensoring with $\mathcal L$ we thus get an
exact sequence
\begin{equation}\label{defT}
0 \to \mathcal \pi^*\mathcal K_k \otimes \mathcal L \to \pi^*
S^{k-1}\Omega_Y \otimes \Omega_X \otimes \mathcal L \to \mathcal
T_k \to 0\,
\end{equation}
where
$$\mathcal T_k := \frac{\pi^* S^{k-1}\Omega_Y \otimes \Omega_X
\otimes \mathcal L}{\pi^* \mathcal K_k \otimes \mathcal L} \ $$ is
a subbundle of $S^k \Omega_X \otimes \mathcal L$. Arguing as in
\cite[pp. 559--560]{LMP} this allows us to obtain the following
exact diagram:

$$\begin{matrix}  &  & 0 &  & 0 &  & 0 & & \\
 &  & \downarrow &  & \downarrow &  & \downarrow &  &\\
 0 & \to & \mathcal Q_{k-1}^{\vee} & \to & \mathcal P^{k-1}_{X}(\mathcal L)^{\vee} & \to &
\mathcal E_{k-1} & \to & 0\\
 &  & \downarrow &  & \downarrow &  & \downarrow & & \\
 0 & \to & \mathcal Q_k^{\vee} & \to & \mathcal P^k_X(\mathcal L)^{\vee} & \to & \mathcal E_k & \to & 0\\
 &  & \psi \downarrow &  & \downarrow &  & \beta \downarrow &  & \\
 0 & \to & \mathcal M_k^{\vee} & \to & (S^k\Omega_X \otimes \mathcal L)^{\vee} & \to &
 \mathcal T_k^{\vee} & \to & 0\\
 &  & \downarrow  &  & \downarrow &  & \downarrow & &  \\
 &  & 0 &  &  0  & & 0 & &.\\
\end{matrix}$$
In particular, the third vertical sequence above,
\begin{equation*}
0 \to \mathcal E_{k-1} \to \mathcal E_k \to \mathcal T_k^{\vee}
\to 0\, ,
\end{equation*}
shows by induction on $k$ that $\mathcal E_k$ is a vector bundle
(note that $\mathcal E_0=\mathcal L^{-1}$ and
$\mathcal E_1=\mathcal P_X^1(\mathcal L)^\vee$).
\end{proof}

\begin{thm}\label{thm-2}
Let $X\subset \mathbb P^N$ be a $n$-dimensional scroll over a smooth
$m$-dimensional variety $Y$, with hyperplane bundle
$\mathcal L=\mathcal O_{P^N}(1)|_X$. Let  $\pi\colon X\to Y$ denote the projection, and set $\mathcal V :=\pi_*\mathcal L$.
For given $k\ge 1$, set  $r_k=(n-m)\binom{m+k-1}{m} + \binom{m+k}m $.
Assume that
the generic rank of $j_k$ is equal to $r_k$ for all $k\ge 1$ such
that $r_k \le N+1$. Assume the $k$th inflectional locus $\Phi_k$
of $X$ has codimension $\ell:=N+2-r_k$ or is empty.
Then $\Phi_k$ has a natural structure as a Cohen--Macaulay
scheme and its class is given as
\[[\Phi_k]=[\prod_{i=0}^{k-1} \pi^*c(S^{i}T_Y \otimes
\mathcal V^{\vee})^{-1}\  c(\pi^*S^kT_Y \otimes \mathcal L^{-1})^{-1}]_\ell.\]
\end{thm}

\begin{proof}
The $k$th inflectional locus of $X$, i.\,e., the locus $\Phi_k
\colon = \Phi_k(X)$ where $j_{k}\colon V_X\to \mathcal
P_X^k(\mathcal L)$ drops rank, is the same as the locus where the
map $j_{k}^\vee\colon \mathcal E_k \to V_X^\vee$ drops rank.
Locally the map is given by a $(N+1) \times r_k$ matrix. Hence the
expected codimension of $\Phi_k$ is $\ell:=N+1-(r_k-1)$. Then the
first part of the statement in the theorem and the fact that
$[\Phi_k] = [c(\mathcal E_k)^{-1}]_{\ell}$ follows from Porteous'
formula \cite[Example 14, p.~255]{Fu}, observing that
$c(V_X^\vee)c(\mathcal E_k)^{-1}=c(\mathcal E_k)^{-1}$.

To make this expression explicit, we then compute the total Chern
class $c(\mathcal E_k)$ by recursion. Since $\mathcal E_0 =
\mathcal L^{-1}$, we get
\begin{equation}\label{provis}
c(\mathcal E_k) = \prod_{i=1}^k c(\mathcal T_i^{\vee}) \
c(\mathcal L^{-1}).
\end{equation}
To compute the total Chern class of $\mathcal T_i^{\vee}$,
dualizing (\ref{defT}) we get
$$
0 \to \mathcal T_i^{\vee} \to \pi^*S^{i-1}T_{Y}^{i}\otimes
T_{X}\otimes \mathcal L^{-1} \to \pi^* \mathcal K_i^{\vee} \otimes
\mathcal L^{-1} \to 0,
$$
where $T_Y$ and $T_X$ are the tangent bundles of $Y$ and $X$.
Hence
\begin{equation}\label{Tdual}
c(\mathcal T_i^{\vee}) = c(\pi^*S^{i-1}T_{Y}^{i}\otimes
T_{X}\otimes \mathcal L^{-1})\ c(\pi^* \mathcal K_i^{\vee}
\otimes \mathcal L^{-1})^{-1}.
\end{equation}
Note that $X=\mathbb P(\mathcal V)$, since $\mathcal V = \pi_*\mathcal
L$, so that $\mathcal L$ is the tautological line bundle on $X$.
Recall the relative tangent and Euler exact sequences
\begin{equation}\label{reltang}
0 \to T_{X/Y} \to T_X \to \pi^* T_Y \to 0,
\end{equation}
\begin{equation}\label{relEuler}
0 \to \mathcal O_X \to \pi^*\mathcal V^{\vee} \otimes \mathcal L
\to T_{X/Y} \to 0,
\end{equation}
Now, twisting (\ref{reltang}) by $\pi^*S^{i-1}T_Y \otimes \mathcal
L^{-1}$, we get the exact sequence:
\setlength\arraycolsep{2pt}
\begin{eqnarray*} 0 \to \pi^*S^{i-1}T_Y \otimes T_{X/Y} \otimes \mathcal L^{-1}
 \to  \pi^*S^{i-1}T_Y \otimes T_X \otimes \mathcal L^{-1}\\
\to \pi^*(S^{i-1}T_Y \otimes T_Y) \otimes \mathcal L^{-1} \to 0.
\end{eqnarray*}
On the other hand, twisting by $\mathcal L^{-1}$ the dual of
(\ref{defK}) lifted to $X$ via $\pi^*$, we obtain a second exact
sequence:
$$0 \to \pi^*S^iT_Y \otimes \mathcal L^{-1} \to
\pi^*(S^{i-1}T_Y \otimes T_Y) \otimes \mathcal L^{-1} \to \pi^*
\mathcal K_i^{\vee} \otimes \mathcal L^{-1} \to 0.$$ Using both,
(\ref{Tdual}) can be rewritten as
$$c(\mathcal T_i^{\vee}) = c(\pi^*S^{i-1}T_Y \otimes T_{X/Y} \otimes \mathcal L^{-1})
c(\pi^*S^iT_Y \otimes \mathcal L^{-1}).$$
Now, twisting
(\ref{relEuler}) by $\pi^*S^{i-1}T_Y \otimes \mathcal L^{-1}$ we
can replace the first total Chern class appearing on the right
hand of the above equality, obtaining
$$c(\mathcal T_i^{\vee}) = \pi^*c(S^{i-1}T_Y \otimes \mathcal V^{\vee})
c(\pi^*S^{i-1}T_Y \otimes \mathcal L^{-1})^{-1}
c(\pi^*S^iT_Y \otimes \mathcal L^{-1}).
$$
Note that the last two terms cancel in the product $\prod_{i=1}^k
c(\mathcal T_i^{\vee})$, except the second for $i=1$, which is
$c(\mathcal L^{-1})^{-1}$, and the third for $i=k$.
Thus, because of these cancellations, (\ref{provis}) gives
\begin{equation}\label{explicit}
c(\mathcal E_k) = \prod_{i=1}^k \pi^*c(S^{i-1}T_Y \otimes
\mathcal V^{\vee})  c(\pi^*S^kT_Y \otimes \mathcal L^{-1}),
\end{equation}
which shows what we wanted.
\end{proof}

Note that we can also write, using Segre classes,
\[c(\mathcal E_k)^{-1}=s(\mathcal E_k^\vee)= \prod_{i=1}^k \pi^*s(S^{i-1}\Omega_Y \otimes
\mathcal V)  s(\pi^*S^k\Omega_Y \otimes \mathcal L).\]

\medskip
As to the assumptions in Theorem \ref{thm-2}, note that the
requirement that $\Phi_k$ has the expected codimension $\ell$ is
independent of the condition $s_k=r_k$, as the following example
shows.
\medskip

\noindent{\it{Example 3.}} Let $\mathcal{V}$ a vector bundle of
rank $2$ on $\mathbb{P}^2$ fitting into an exact sequence
$$0 \to \mathcal{O}_{\mathbb{P}^2}(1)^{\oplus 2} \to
T_{\mathbb{P}^2}^{\oplus 2} \to \mathcal{V} \to 0.$$ Then
$\mathcal{V}$ is very ample and $c_1(\mathcal{V})=4$,
$c_2(\mathcal{V})=6$. Letting $X:=\mathbb{P}(\mathcal{V})$, the
tautological line bundle $\mathcal{L}$ embeds $X$ in
$\mathbb{P}^9$ as a scroll of degree $10$. In fact $X$ can be
regarded as a general threefold section of $\mathbb{P}^2 \times
\mathbb{P}^3 \subset \mathbb{P}^{11}$, embedded via the Segre
embedding, with a linear space \cite[Section 3]{LN}. Sometimes one
refers to $X$ as a Bordiga scroll, since its general hyperplane
section $S \in |\mathcal{L}|$ is a Bordiga surface, namely a
rational surface of sectional genus $3$. According to
\cite[Theorem 3]{Ma}, $X$ is also isomorphic to $\mathbb{P}^3$
blown-up along a twisted cubic $C$ and, by \cite[Lemma 1.4]{LN},
$\mathcal{L}=\sigma^*\mathcal{O}_{\mathbb{P}^3}(3) \otimes
\mathcal{O}_X(-E)$, where $\sigma$ is the blowing-up and $E$ is
the exceptional divisor. Thus $V=H^0(X,\mathcal{L}) \cong
H^0(\mathbb{P}^3,\mathcal{J}_C(3))$, where $\mathcal{J}_C$ is the
ideal sheaf of $C \subset \mathbb{P}^3$. We study $\Phi_2(X)$
relying on this description. Using coordinates $(1:x:y:z)$ on
$\mathbb{P}^3$ let $C$ be defined by $(1:t:t^2:t^3)$. Consider the
three quadrics $Q_1, Q_2, Q_3$ containing $C$, defined by $x^2-y$,
$xy-z$, $xz-y^2$ respectively. Taking into account the syzygies
$xQ_2=yQ_1-Q_3$ and $zQ_1 = -yQ_2+xQ_3$, we can easily find a
basis for $H^0(\mathbb{P}^3,\mathcal{J}_C(3))$. Moreover, as
$\sigma$ is an isomorphism outside $E$ we can use the affine
coordinates $(x,y,z)$ around any point of $\mathbb{P}^3 \setminus
C$ as local coordinates around the corresponding point $p \in X
\setminus E$. Then $X \setminus E$ can be represented in
$\mathbb{P}^9$, near $p$, by the following homogeneous coordinates
$$(Q_1: xQ_1: yQ_1: Q_2: yQ_2: zQ_2: Q_3: xQ_3: yQ_3: zQ_3).$$
Now computing the rank of $j_{2,p}$ (e.\,g., by using \emph{Maple}), we
obtain that $s_2=9=r_2$. Actually, computing the minors of order
$9$ one can easily find that the locus they define is $y=x^2,
z=x^3$. This means that $\text{rk}\,j_{2,p}=9$ for every $p \in X
\setminus E$. Hence $\Phi_2 \subseteq E$. To look at $E$, note
that the plane $\Pi \subset \mathbb{P}^3$ defined by $x=0$ is
transverse to $C$ at the point $o=(0,0,0)$. Taking $(y,w)$ as
local coordinates on $\widetilde{\Pi}:=\sigma^{-1}(\Pi)$, the
local equations of $\sigma$ restricted to $\widetilde{\Pi}$ are
$\sigma(0,y,w)=(0,y,z=yw)$. Therefore the homogeneous coordinates
in $\mathbb{P}^9$ of points $p \in X$, near the fiber $f_o \subset
E$ of $\sigma|_E$, are given by the following polynomials: $x^2-y,
x^3-xy, x^2y-y^2, xy-yw, xy^2-y^2w, xy^2w-y^2w^2, xyw-y^2,
x^2yw-xy^2, xy^2w-y^3, xy^2w^2-y^3w$, in terms of the local
parameters $x,y,w$ (simply one has to replace $z$ with $yw$ in the
previous polynomials). Now redoing the computation (with \emph{Maple}) we
see that the  rank of $ j_{2,p}$ is $9$, for $p$
in a neighborhood $U$ of $f_o$ in $X$. But all $9\times 9$-minors
contain the factor $y$ with some multiplicity. Since $y=0$ is the
local equation of $E$ in $U$, this means that
$(\Phi_2)_{\text{red}} = E$. In particular $\Phi_2$ is a divisor.
We conclude that $s_2=r_2$, but that $\Phi_2$ does not have the expected
codimension $\ell=2$.

\medskip
The assumption that the inflectional locus has codimension $\ell$
implies that $\ell$ is in the range $1\le \ell \le n$, which means
\[1\le N+2-r_k\le n.\]
This implies that
$k$ is the maximal integer such that $r_k\le N+1$ holds, as
we will now explain. First of all, $\ell\ge 1$ is equivalent
to $r_k\le N+1$, and
$\ell\le n$ is equivalent to $r_k\ge N+2-n$. If this holds, then since
\[\textstyle r_{k+1}=r_k+\binom{m+k}{m-1}+(n-m)\binom{m-1+k}{m-1}\]
we get
\[\textstyle r_{k+1}\ge N+2-n +
\binom{m+k}{m-1}+(n-m)\binom{m-1+k}{m-1}.\]
Hence it suffices to show that
\[\textstyle -n +
\binom{m+k}{m-1}+(n-m)\binom{m-1+k}{m-1}\ge 0.\]
The left hand side can be written as
\[\textstyle n(\binom{m-1+k}{m-1}-1)-k\binom{m-1+k}{m-2}\]
which is the same as
\[\textstyle k\big((n-1)+(\frac{n}2-1)(k+1)+\cdots + (\frac{n}{n-m}-1)\binom{m-2+k}{m-2}\big),\]
and this is clearly non-negative.
Hence $r_{k+1}\ge N+2$.

In the case that $Y$ is a curve (see \cite{LMP}), so that $m=1$,
then $r_k=kn+1$, and in this case, if $k$ is the largest integer
such that $r_k\le N+1$, then conversely $1\le \ell \le n$ holds.
When $m\ge 2$, this last implication is no longer true. However,
if we take $k$ to be the largest integer such that $r_{k-1}-1+n\le
N$, then the bound on $\ell \le n$ holds (recall that $r_{k-1}-1$
is the dimension of the $(k-1)$th osculating space to $X$ at a
general point). To see this, assume $r_k-1+n\ge N+1$. Then
$r_k\ge N+2-n$, so that
\[\ell=N+2-r_k \le N+2 - N-2+n=n.\]

The scrolls for which Theorem \ref{thm-2} applies are thus those
such that the general $k$th osculating spaces are of dimension
$r_k-1$ and where $r_k$ satisfies
\begin{equation}\label{range}
r_k-1\le N\le r_k+n-2.
\end{equation}
For instance, the appropriate range for studying $\Phi_2$ is
\begin{equation}\label{range,k=2}\textstyle
(n-m)(m+1)+\binom{m+2}2-1\leq N\leq (n-m)(m+1)+\binom{m+2}2+n-2.
\end{equation}
For $m=2$ this gives
\[3n-1\leq N\leq 4n-2,\]
and for $m=3$,
\[4n-3 \leq N\leq  5n-4.\]

In the special cases $m=2$, $n=3$ and $m=3$, $n=4$ that we will
consider in Sections \ref{k=2} and \ref{k=2,m=3}, this means that $8 \leq N \leq 10$ and
$13\le N \le 16$ respectively.

Note however that (\ref{range}) in no way implies that $s_k=r_k$. For instance, for $(\mathbb{P}^1)^4 \subset \mathbb{P}^{15}$
condition (\ref{range,k=2}) is satisfied but, as observed in
Example 2, we have $s_2=11<14=r_2$.

Sometimes in the following we say that $X \subset \mathbb{P}^N$
\emph{satisfies our general assumptions} to mean that $s_k=r_k$
and $\Phi_k$
has the expected codimension. From Section \ref{k=2,divisor} on, we in addition assume $k=2$.

\section{The case $Y$ is an abelian variety}
In this section we assume that $X\subset \mathbb P^N$ is a scroll over an abelian variety $Y$.
We use the same notations as in the previous section: $\pi\colon X\to Y$ is the scroll projection, $\mathcal L$ is the hyperplane bundle on $X$, and $\mathcal V=\pi_* \mathcal L$.

In this case, we have $T_{Y}=\mathcal
O_{Y}^{\oplus m}$, so $S^jT_{Y}=\mathcal O_{Y}^{\oplus
\binom{m-1+j}{m-1}}$ and $\pi^*S^jT_Y=\mathcal
O_X^{\oplus\binom{m-1+j}{m-1}}$. Set $\mu:=\binom{m-1+k}{m-1}.$
Then
$$ \textstyle c(\pi^* S^kT_Y \otimes \mathcal L^{-1})=
c(\pi^* \mathcal O_Y^{\oplus \mu} \otimes \mathcal L^{-1}) =
c\big({(\mathcal L^{-1})}^{\oplus \mu}\big) = (1-L)^{\mu},$$ where
we set $L:=c_1(\mathcal L)$. Similarly,
$$\prod_{i=0}^{k-1} c(S^{i}T_Y\otimes \mathcal V^{\vee}) =
\prod_{i=0}^{k-1} c(\mathcal V^{\vee})^{\binom{m-1+i}{m-1}}=
c(\mathcal V^{\vee})^{\nu},$$ where $\nu :=
\sum_{i=0}^{k-1}\binom{m-1+i}{m-1}=\binom{m-1+k}m$. Therefore we
get from (\ref{explicit})
$$c(\mathcal E_{k})= (1-L)^{\mu}  \pi^* c(\mathcal V^{\vee})^{\nu}.$$
Now we have to compute the inverse
\begin{equation*}
c(\mathcal E_{k})^{-1}=(1-L )^{- \mu} \pi^*c(\mathcal V^{\vee})^{-
\nu}.
\end{equation*}
By the general binomial formula,
$$\textstyle(1-L )^{- \mu} = \sum_{j=0}^{\infty} \binom{j+ \mu -1}{\mu-1} L^j.$$

Set $V_i = \pi^*c_i(\mathcal V)$. We have
\setlength\arraycolsep{2pt}
\begin{eqnarray*}\textstyle
\pi^* c(\mathcal V^{\vee})^{- \nu} \nonumber &=&
(1-V_1+V_2-V_3+\cdots)^{- \nu} \nonumber \\
&=&\textstyle \sum_{h=0}^\infty \binom{h+ \nu -1}{\nu -1}
(V_1-V_2+\cdots +(-1)^{m+1}V_m)^h.
\end{eqnarray*}
All this gives
\begin{equation}\label{inverse-fin}\textstyle
c(\mathcal E_{k})^{-1}=\sum_{j=0}^n \binom{j+ \mu -1}{\mu -1} L^j
\sum_{h=0}^m \binom{h+ \nu -1}{\nu -1} (V_1-V_2+\cdots
+(-1)^{m+1}V_m)^h .
\end{equation}
\medskip

Now let us put $m=2$. Then we get the following result.

\begin{prp} \label{Theorem-abelian}
Let $X \subset \mathbb P^N$ be a $n$-dimensional
scroll over an abelian surface $Y$, with projection $\pi: X \to
Y$, and let $\mathcal V=\pi_*\mathcal L$, where $\mathcal L$ is
the hyperplane bundle. Suppose that $(\ref{range})$ is satisfied,
that $j_{k}$ has maximal general rank, and that the inflectional
locus $\Phi_k$ has the expected codimension $\ell$ or is empty. Then the
cohomology class of $\Phi_k$ is
\setlength\arraycolsep{2pt}
\begin{eqnarray*}
[\Phi_k]& =&\textstyle  \binom{\ell+k}k L^{\ell} + \nu \binom{\ell-1+k}k V_1
L^{\ell - 1}
 + \textstyle \binom{\ell-2+k}k \left[ \binom{\nu +1}2
V_1^2 - \nu V_2 \right] L^{\ell -2},
\end{eqnarray*}
where $L=c_1(\mathcal L)$ and $V_i = \pi^* c_i(\mathcal V)$. In
particular,
\setlength\arraycolsep{2pt}
\begin{eqnarray*}
\deg \Phi_k  &=&\textstyle \left[\binom{\ell +k}k+ \binom{k+1}2 \binom{\ell -2
+ k}k \right] d\\
&&\textstyle + \left[\binom{k+1}2 \binom{\ell -1+k}k +
\binom{\binom{k+1}2}2 \binom{\ell -2+k}k \right] (2g-2),
\end{eqnarray*}
where
$d$ and $g$ are the degree and the sectional genus of $X$.
\end{prp}
\begin{proof} The formula for $[\Phi_k]$ follows from Theorem \ref{thm-2} and
(\ref{inverse-fin}). To get the degree it is enough to dot with
$L^{n-\ell}$ and to recall the following facts: $d= L^n =
L^{n-2}(V_1^2-V_2)$ and $ L^{n-1} \pi^*D = c_1(\mathcal V) D$ for
any divisor $D$ on $Y$, by the Chern--Wu relation. Moreover, for
the general surface section $S$ of $X$, the pair $(Y, \det
\mathcal V)$ is the adjunction theoretic reduction of $(S,\mathcal
L_S)$ via the map $\pi|_S:S \to Y$. Then $g = g(X,\mathcal L) =
g(S, \mathcal L_S) = g(Y, \det \mathcal V)$, hence $c_1(\mathcal
V)^2 = 2g-2$ by the genus formula.
\end{proof}

In particular, suppose that $\Phi_k$ has codimension $\ell=n$,
which implies $N=n(1+ \nu)-1-\binom{k}2$. Then $X$ has finitely
many flexes, their number being given by
\setlength\arraycolsep{2pt}
\begin{eqnarray}\label{abelian}\deg \Phi_k &=&\textstyle \left[\binom{n +k}k + \binom{k+1}2 \binom{n-2+k}k\right]d
\nonumber \\
&&\textstyle + \left[\binom{k+1}2 \binom{n -1+k}k + \binom{\binom{k+1}2}2
\binom{n -2+k}k\right](2g-2).
\end{eqnarray}
 For instance, for $\ell=n=3$ and
$k=2$, we get $\deg \Phi_2 = 19 d + 27 (2g-2)$.  Note that $d$ and
$g$ are large; actually $\det \mathcal V$ is a very ample line
bundle on $Y$, since $\mathcal V=\pi_*\mathcal L$ is a very ample
vector bundle. So, in a sense, the above formula shows that any
$n$-dimensional scroll $X$ over an abelian surface is very, very
far from being uninflected.
\medskip

{\it{Example 4}}. Let $A$ be an abelian surface with Picard number
$1$, and let $\mathcal{H}$ be a line bundle on $A$ representing a
polarization of type $(1,p)$ with $p=(k+1)^2+2$, $k \geq 2$. Then
$\mathcal{H}$ is very ample by Reider's theorem. Moreover,
$H^0(A,\mathcal{H})$ embeds $A$ in $\mathbb{P}^N$ with $N=p-1$ as
a surface of degree $2p$. Let $\epsilon \in A$ be a nontrivial
element of order $2$, let $Y$ be the quotient $A/\epsilon$, and let
$q:A \to Y$ be the quotient map, which is an isogeny of degree
$2$. Set $\mathcal{V}:=q_*\mathcal{H}$. Then $\mathcal{V}$ is a
rank-2 vector bundle on $Y$. As $q^*\mathcal{V}= \mathcal{H}
\oplus t_{\epsilon}^*\mathcal{H}$, where $t_{\epsilon}:A \to A$ is
the translation by $\epsilon$ \cite[(1)]{CH1}, we get
$$\big(c_1(\mathcal{V})^2,c_2(\mathcal{V})\big) =
(4p,p)$$ by the functoriality of the Chern classes. Now let
$X:=\mathbb{P}(\mathcal{V})$ and let $\mathcal{L}$ be the
tautological line bundle. Note that $p > 9$. Then \cite[Theorem
3.3]{CH1} implies the following facts: $\mathcal{L}$ is very
ample, $H^0(X,\mathcal{L})$ embeds $X$ in $\mathbb{P}^N$, and
identifying this space with that containing $A$ via the
isomorphisms
$$H^0(A,\mathcal{H}) \cong H^0(Y, \mathcal{V}) \cong H^0(X, \mathcal{L}),$$
the image is the $3$-dimensional scroll over $Y$ generated by the
secant lines joining $x$ and $x+\epsilon$, as $x$ varies on $A$.
Note that its degree is
$d=c_1(\mathcal{V})^2-c_2(\mathcal{V})=3p$. Moreover, condition
(\ref{range}) is satisfied with equality on the right, because
$r_k+n-2 = \binom{k+1}2 + \binom{k+2}2 + 1 = (k+1)^2+1 = N$. Thus
for this abelian scroll our formula allows us to express $\deg
\Phi_k$ as a degree $5$ polynomial in $k$, provided that $\Phi_k$
is a finite set. The explicit expression is the following: $\deg
\Phi_k=\frac{p}{2} (k^5 + 5k^4 + 13k^3 + 19k^2 + 16k +6)$. In
particular, for $k=2$ we get $\deg \Phi_2 =1815$. This example can
be generalized to scrolls of any dimension $n \geq 3$ over an
abelian surface by taking $p=n \big(\binom{k+1}2 + 1\big) -
\binom{k}2$ and $Y: = A/G$, where $G \subset A$ is the subgroup
generated by a nontrivial element $\rho$ of order $n-1$. In this
case $\mathcal{V}=q_*\mathcal{H}$ is a vector bundle of rank $n-1$
and $q^*\mathcal{V}= \mathcal{H} \oplus t_{\rho}^*\mathcal{H}
\oplus \dots \oplus t_{(n-2)\rho}^*\mathcal{H}$. As before, let
$X:= \mathbb{P}(\mathcal{V})$, so that $X$ is an $n$-dimensional
scroll over $Y$, and let $\mathcal{L}$ be the tautological line
bundle again. Since $A$ has Picard number $1$ and the inequality
$\mathcal{H}^2 \geq 4(2n-1)+6$ is satisfied for every $n \geq 3$
and $k \geq 2$, we know from \cite[Theorem 1.1]{T} that
$\mathcal{H}$ is $(2n-1)$-very ample on $A$. Then
\cite[Proposition 2.2]{CH2} tells us that $\mathcal{L}$ is very
ample, $H^0(X,\mathcal{L})$ embeds $X$ in $\mathbb{P}^N$, and
identifying this space with that containing $A$ as before, the
image is the $n$-dimensional scroll over $Y$ generated by the
$(n-2)$-dimensional linear spaces $\langle x, x+\rho, \dots , x +
(n-1)\rho \rangle$, as $x$ varies on $A$. Condition (\ref{range})
is satisfied and $\ell=n$ again. Then our formula for $\deg
\Phi_k$ applies also in this case, provided that $\Phi_k$ is a
finite set.
\medskip

For $m=3$ we get the following result.

\begin{prp} \label{Theorem-abelian3}
Let $X \subset \mathbb P^N$ be a $n$-dimensional
scroll over an abelian threefold $Y$, with projection $\pi: X \to
Y$, and let $\mathcal V=\pi_*\mathcal L$, where $\mathcal L$ is
the hyperplane bundle. Suppose that $(\ref{range})$ is satisfied,
that $j_{k}$ has maximal general rank, and that the inflectional
locus $\Phi_k$ has  the expected codimension $\ell$. Set
$L=c_1(\mathcal L)$ and $V_i = \pi^*c_i(\mathcal V)$. Then the
cohomology class of $\Phi_k$ is
\setlength\arraycolsep{2pt}
\begin{eqnarray*}
[\Phi_k]& =& \textstyle \binom{\ell+ \mu -1}{\mu -1} L^{\ell} + \nu
\binom{\ell+ \mu -2}{\mu -1} V_1 L^{\ell - 1} \\ && \textstyle + \binom{\ell+
\mu -3}{\mu -1} \big( \binom{\nu +1}2 V_1^2 - \nu V_2 \big)
L^{\ell -2}\\
&&\textstyle +\binom{\ell + \mu -4}{\mu -1}\big(\binom{\nu +2}{3}V_1^3-
\nu (\nu +1)V_1V_2+ \nu V_3\big)L^{\ell -3},
\end{eqnarray*}
 where $\mu
=\binom{k+2}2$ and $\nu =\binom{k+2}3$.
\end{prp}

\section{The case $k=2$: when $\Phi_2$ is a divisor}
\label{k=2,divisor}

From now on we consider the case $k=2$. In this Section we
study the inflectional locus $\Phi_2$, assuming it is a divisor. First
of all we prove the following:
\begin{thm} \label{thm-divisor} Let
$X\subset \mathbb P^N$ be a $n$-dimensional scroll over a smooth
$m$-dimensional variety $Y$, with projection $\pi$ and let
$\mathcal{V}=\pi_*\mathcal{L}$, where $\mathcal L$ is the
hyperplane bundle. Suppose that the conditions on the range of $N$ are
satisfied, that $j_2$ has maximal general rank and that the
inflectional locus $\Phi_2$ is a divisor. Then its class is given
by
\begin{equation*}\textstyle [\Phi_2] = \pi^*\big((n+2)K_Y + (m+1) c_1(\mathcal{V})\big) +
\binom{m+1}2 L,
\end{equation*}
 where $L = c_1(\mathcal{L})$.
\end{thm}
\begin{proof}
Letting $k=2$, (\ref{explicit}) reduces to the product of three
terms
\begin{equation} \label{tag 8}
 c(\mathcal E_2)= \pi^* c(\mathcal V^{\vee})\ \pi^* c(\mathcal V^{\vee}
\otimes T_Y)\ c(\pi^*S^2 T_Y \otimes \mathcal L^{-1}).
\end{equation}
According to Theorem \ref{thm-2}, to determine the class of
$\Phi_2$ we have to take the inverse of the expression in
(\ref{tag 8}). If $\Phi_2$ is a divisor, i.\,e., $\ell = 1$, we need
only the degree 1 term. In other words,
$$[\Phi_2] = \alpha_1 + \beta_1 + \gamma_1,$$
where $\alpha_1, \beta_1$ and $\gamma_1$ are the terms of degree 1
appearing in the expressions of $c(\pi^* \mathcal
V^{\vee})^{-1}$, $c(\pi^* \mathcal V^{\vee} \otimes
T_Y)^{-1}$, and $c(\pi^* S^2T_Y \otimes
\mathcal{L}^{-1})^{-1}$, respectively. Thus
$$\alpha_1 = -c_1(\pi^* \mathcal V^{\vee}) = \pi^* c_1(\mathcal V).$$
Similarly, $\beta_1 = -c_1(\pi^* \mathcal V^{\vee} \otimes T_Y)$.
By confining the use of the splitting principle to $\mathcal{V}$,
we can easily see that
$$\beta_1 = (n-m+1) \pi^*K_Y + m \pi^* c_1(\mathcal V).$$
Finally, $\gamma_1 = -c_1(\pi^* S^2T_Y \otimes \mathcal{L}^{-1})$.
Noting that the rank of $S^2T_Y$ is $\binom{m+1}2$ and
$c_1(S^2T_Y)=(m+1)c_1(T_Y)$ we get
$$\textstyle \gamma_1 = (m+1)\pi^*K_Y + \binom{m+1}2 L.$$
Adding the three terms, we thus get the expression in the
statement.
\end{proof}

\medskip
\noindent{\it{Example 5}.} Consider the scroll over $\mathbb{P}^m$
given by $X= \mathbb{P}(\mathcal{V})$, where $\mathcal{V}=
\mathcal{O}_{\mathbb{P}^m}(1) \oplus
\mathcal{O}_{\mathbb{P}^m}(2)$. Here $n=m+1$ and the tautological
line bundle $\mathcal{L}$ embeds $X \subset \mathbb{P}^N$ where $N
= \binom{m+2}m + m$. Note that, according to (\ref{range,k=2}), we
are in the appropriate range for studying $\Phi_2$ with $\ell=1$,
so that $\Phi_2$ is a divisor. Let $Y_0 \subset X$ be the section
corresponding to the summand $\mathcal{O}_{\mathbb{P}^m}(1)$.
Equivalently we can look at $Y_0$ as the tautological section of
the vector bundle $\mathcal{V} \otimes
\mathcal{O}_{\mathbb{P}^m}(-2)$; hence $Y_0$ is negative and $L =
Y_0 + 2\pi^*\mathcal{O}_{\mathbb{P}^m}(1)$. By applying
Proposition \ref{thm-divisor} we thus get
$$\textstyle \Phi_2(X) = \binom{m+1}2 Y_0.$$
This allows us to point out another relevant discrepancy with
respect to the case of scrolls over curves. Let $m=1$. Then our
$X$ is the cubic surface scroll $X \subset \mathbb{P}^4$ and $Y_0$
is its directrix line. Recall that $r_2=5$ for $m=1$.
Parameterizing $X$, locally around $Y_0$, by
$$(1: u: v: vu: vu^2)$$
we see that the inflectional locus $\Phi_2$ is defined by $v=0$,
hence it consists of $Y_0$, whose codimension is the expected one,
and $\text{rk}\, j_{2,x}=4$ for every $x \in \Phi_2$, i.\,e., one
less than the generic rank. By applying \cite[Corollary 1]{LMP} to
this case we get $[\Phi_2]= L-2F$, where $L$ and $F$ are the
classes of a hyperplane and a fiber respectively. And this agrees
with the fact that $Y_0 \in |L-2F|$. Now let $m=2$ and note that
here $r_2=9$. Parameterizing $X$, locally around $Y_0$, by
$$(1: u_1: u_2: v: vu_1: vu_2: vu_1^2: vu_1u_2: vu_2^2)$$
a direct check shows that $s_2=9$ and $\Phi_2$ is defined by
$v^3=0$. So $\Phi_2$ is non reduced. In fact $\Phi_2(X)=3Y_0$,
according to the above formula. Set-theoretically, $\Phi_2= Y_0$
and is of the expected codimension, but $\text{rk}\,j_{2,x}=6=r_2-3$
for every $x \in \Phi_2$. In other words, the rank of $j_{2}$
drops by $3$ on the inflectional locus. From the geometric point
of view this fact can be rephrased as follows. Take any point $x
\in Y_0$, let $\{\lambda_t\}$ be the pencil of lines in $Y_0$ passing
through $x$, and set $S_t=\mathbb{P}(\mathcal{V}_{\lambda_t})$. Then
$\text{Osc}^2_x(X)= \big\langle \bigcup_t \text{Osc}^2_x(S_t)
\big\rangle$. Note that, for every $t$, $\text{Osc}^2_x(S_t)$ is a
$\mathbb{P}^3$ containing the fiber of $X$ through $x$. This says
that $\text{Osc}^2_x(X)$ is a $\mathbb{P}^5$. Moreover, we see
that $Y_0$, which is the union of the directrix lines constituting
the $\Phi_2(S_t)$ as $\lambda_t$ varies, is equal to
$\Phi_2(X)_{\text{red}}$.
\medskip

\noindent{\it{Remark}}. Suppose $\Phi_2$ is a divisor. Theorem
\ref{thm-divisor} implies that
$$\textstyle (\Phi_2)_f \in
\big|\mathcal{O}_{\mathbb{P}^{n-m}}(\binom{m+1}2)\big|$$ for any
fiber $f$ of $X$. This prevents $\Phi_2$ from being empty, since
by definition $\Phi_2$ is either effective or empty. Therefore $X$ is
not uninflected. Moreover, since $f$ is a linear
$\mathbb{P}^{n-m}$, for every line $\lambda \subset f$ the condition
$\Phi_2 \,  \lambda = \binom{m+1}2$ implies that
$$\textstyle \deg \Phi_2 \geq \binom{m+1}2.$$
Otherwise, every
fiber $f$ would be contained in $\Phi_2$, implying that
$\Phi_2(X)=X$, a contradiction. Note that for the scroll
considered in the above example we have in fact $\deg \Phi_2 =
\binom{m+1}2$.

\medskip
Dotting the expression of $\Phi_2$ provided by Theorem
\ref{thm-divisor} with $L^{n-1}$ we get
\begin{equation} \label{deg Phi2}\textstyle
\deg \Phi_2 =
\pi^*\big((n+2)K_Y+(m+1)c_1(\mathcal{V})\big)L^{n-1}+ \binom{m+1}2
d, \end{equation} where $d$ is the degree of $X$. For instance,
let $m=2$, so that $\mathcal{V}$ has rank $n-1$. As we observed in
the proof of Proposition \ref{Theorem-abelian}, the sectional
genus of $X$ is $g=g(Y, \det \mathcal{V})$. Then, taking into
account the Chern--Wu relation and genus formula, we get the
following two equivalent expressions for the degree:
\setlength\arraycolsep{2pt}
\begin{eqnarray}\label{two-expressions}
\deg \Phi_2 &=& (4-n)d+ (n+2)(2g-2) -
(n-1)c_2(\mathcal{V})\nonumber\\
&= &(4-n)c_1(\mathcal{V})^2 + (n+2)(2g-2) - 3c_2(\mathcal{V}) .
\end{eqnarray}
\medskip
From (\ref{deg Phi2}) we immediately get the following
consequence.

\begin{cor}\label{cor7} Let $X$, $\mathcal{L}$, $Y$ and $\mathcal{V}$ be as in
Theorem \ref{thm-divisor}. If the $\mathbb{Q}$-divisor $K_Y +
\frac{m+1}{n+2} \det \mathcal{V}$ is nef, then
$$\textstyle \deg \Phi_2 \geq \binom{m+1}2 d,$$
where $d$ is the degree of $X$.
\end{cor}
Note that $\frac{m+1}{n+2}\leq \frac{m+1}{m+3} = 1 -
\frac{2}{m+3}$. Hence the adjunction bundle appearing in the
statement does not fit with the range of the nef values of adjoint
bundles to an ample vector bundle considered by Ohno in \cite{Oh}.
This prevents us from giving a more detailed statement in the
general case. However, restricting to threefold scrolls over surfaces,
we can produce a precise result in the same vein as
Corollary \ref{cor7}. This will be done in the next section.

\section{The divisor $\Phi_2$ for threefold scrolls over surfaces}
\label{k=m=2,divisor}

Let $X=\mathbb{P}(\mathcal{V})$ be a threefold scroll over a smooth
surface $Y$, satisfying our general assumptions and with $\ell=1$,
which implies that $X \subset \mathbb{P}^8$. The expression
provided by Theorem \ref{thm-divisor} becomes
\begin{equation} \label{divisor}
[\Phi_2] = \pi^*\big(5K_Y + 3c_1(\mathcal{V})\big) + 3L.
\end{equation}
\medskip
We start by recalling the following well-known facts (e.\ g., for A2), see
\cite[Proposition 5.2]{Io}).

\medskip

\noindent{\it{Remark}}. Let $\mathcal{V}$ be an ample vector
bundle of rank $2$ on a smooth projective variety $Y$.

\noindent A1) If $C \subset Y$ is a smooth rational curve, then
$\deg \mathcal V_C \geq 2$, and equality implies that
$\mathcal{V}_C \cong \mathcal{O}_{\mathbb{P}^1}(1)^{\oplus 2}$;
moreover, if $\deg \mathcal V_C = 3$, then $\mathcal{V}_C
\cong \mathcal{O}_{\mathbb{P}^1}(2)\oplus
\mathcal{O}_{\mathbb{P}^1}(1)$.

\noindent A2) Suppose that $\mathcal{V}$ is very ample and let $C
\subset Y$ be a smooth curve of genus $1$; then $\deg \mathcal V_C
\geq 5$.

\medskip
\par
Now, let $Y$ be a smooth surface and let $\mathcal{V}$ be a very
ample vector bundle of rank $2$ on $Y$. We want to investigate the
effectivity of $5K_Y+3 \det \mathcal{V}$. First of all we consider
the adjoint line bundle $A:= K_Y + \det \mathcal{V}$. As usual in
adjunction theory we use the additive notation for the tensor
product of line bundles on $Y$.

\begin{prp} \label{prop-adjunction} The line bundle $A$ is very ample
except in the following cases:
\begin{enumerate}
    \item $(Y, \mathcal{V})=(\mathbb{P}^2,
    \mathcal{O}_{\mathbb{P}^2}(1)^{\oplus 2})$;
    \item $(Y, \mathcal{V})=(\mathbb{P}^2,
    \mathcal{O}_{\mathbb{P}^2}(2)\oplus \mathcal{O}_{\mathbb{P}^2}(1))$;
    \item $(Y, \mathcal{V})=(\mathbb{P}^2,
    T_{\mathbb{P}^2})$, where $T_{\mathbb{P}^2}$ is the tangent bundle;
    \item $(Y, \mathcal{V})= (\mathbb{P}^1\times
    \mathbb{P}^1, \mathcal{O}_{\mathbb{P}^1\times
    \mathbb{P}^1}(1,1)^{\oplus 2})$
    \item $Y$ is a $\mathbb{P}^1$-bundle over a smooth curve $B$
    and $\mathcal{V}_f = \mathcal{O}_{\mathbb{P}^1}(1)^{\oplus 2}$
    for every fiber $f$ of the bundle projection $p:Y \to B$.

\end{enumerate}
\end{prp}
\begin{proof} Consider the pair $(Y, \det \mathcal{V})$. According
to \cite[Theorem 0.1]{SV} and taking into account Remark A1), $A$
is nef unless $(Y, \det \mathcal{V})=(\mathbb{P}^2,
\mathcal{O}_{\mathbb{P}^2}(2))$. This gives exception (1), due to
Remark A1) again and the uniformity of $\mathcal{V}$. Now,
according to \cite[Theorem 0.2]{SV} $A$ is big unless $Y$ is a Del
Pezzo surface with $\det \mathcal{V}=-K_Y$, or $(Y,\det
\mathcal{V})$ is a conic fibration over a smooth curve $B$. In the
former case, $Y$ has to be a minimal surface in view of Remark A1),
hence $(Y, \det \mathcal{V})$ is either $(\mathbb{P}^2,
\mathcal{O}_{\mathbb{P}^2}(3))$, or $(\mathbb{P}^1\times
\mathbb{P}^1, \mathcal{O}_{\mathbb{P}^1\times
\mathbb{P}^1}(2,2))$. This gives exceptions (2), (3) and (4), due
to Remark A1) again and the uniformity of $\mathcal{V}$ (e.\ g.,
see \cite[Theorem 2.2.2]{OSS} and \cite[Lemma 3.6.1]{W}). If $(Y,
\det \mathcal{V})$ is a conic fibration then all fibers of the
projection $p:Y \to B$ are irreducible by Remark A1), hence $Y$ is
in fact a $\mathbb{P}^1$-bundle over $B$ and we get case (5), by
Remark A1) again. So, apart from cases (1)--(5), $A$ is nef and big
and then we can consider the adjunction theoretic reduction of
$(Y, \det \mathcal{V})$ \cite[Theorem 0.3]{SV}, which coincides
with $(Y, \det \mathcal{V})$ itself, once more in view of Remark
A1). Then, according to \cite[Theorem 2.1 and Theorem 1.5]{SV} and
Remark A1), we have that $A$ is very ample, except in the following
cases:
\begin{enumerate}
    \item[(i)] $Y$ is a Del Pezzo surface with $K_Y^2=1$ and
    $\det \mathcal{V}=-3K_Y$;
    \item[(ii)] $Y$ is a Del Pezzo surface with $K_Y^2=2$ and
    $\det \mathcal{V}=-2K_Y$.
    \item[(iii)] $Y$ is an elliptic $\mathbb{P}^1$-bundle of
    invariant $-1$ with $\det \mathcal{V}$ linearly equivalent to $3\sigma$,
    where $\sigma$ is a section of minimal self-intersection.
\end{enumerate}
However, they cannot occur in our setting. Actually, in all these
cases $Y$ contains a smooth curve $C$ of genus $1$ such that $\deg
\mathcal{V}_C < 5$ (take $C \in |-K_Y|$ in cases (i), (ii), and
$C=\sigma$ in case (iii)). This contradicts Remark A2).
\end{proof}

\begin{prp} \label{prop-adjunction2} Suppose that the adjoint line bundle
$A$ is very ample. Then $h^0(K_Y+A)>0$ except in the following
cases:
\begin{enumerate}
    \item $(Y, \det \mathcal{V})=(\mathbb{P}^2,
    \mathcal{O}_{\mathbb{P}^2}(4))$;
    \item $(Y, \det \mathcal{V})=(\mathbb{P}^2,
    \mathcal{O}_{\mathbb{P}^2}(5))$;
    \item $Y$ is a $\mathbb{P}^1$-bundle over a smooth curve $B$,
    and $\mathcal{V}_f = \mathcal{O}_{\mathbb{P}^1}(2) \oplus
    \mathcal{O}_{\mathbb{P}^1}(1)$ for every fiber $f$ of the
    bundle projection $p:Y \to B$.
\end{enumerate}
\end{prp}
\begin{proof} Let $C \in |A|$ be a smooth curve and let
$g':=g(C)$. From the cohomology sequence of
$$0 \to K_Y \to K_Y+A \to (K_Y+A)_C \to 0,$$
recalling that $h^0\big((K_Y+A)_C\big)=h^0(K_C)=g'$, we immediately
get, in view of the Kodaira vanishing theorem,
$$h^0(K_Y+A) = p_g(Y) + g' - q(Y).$$
Note that $g' \geq q(Y)$ by the Lefschetz theorem, with equality
if and only if $(Y,A)$ is either $(\mathbb{P}^2,
\mathcal{O}_{\mathbb{P}^2}(a))$ with $a=1,2$, or a scroll over a
smooth curve $B$ \cite[Corollary 1.5.2]{So}. The former case leads
to (1) and (2) respectively, while in the latter, we have $\deg
\mathcal{V}_f = 3$ for every fiber $f$ of the projection $p:Y \to
B$, since $K_Yf=-2$. Then, taking into account Remark A1), we get
(3). We thus see that $h^0(K_Y+A)\geq 1$ apart from (1)--(3).
\end{proof}

\begin{cor}\label{cor-effectivity}
Suppose that $(Y,\mathcal{V})$ is not any of the exceptions listed
in Proposition \ref{prop-adjunction} and Proposition
\ref{prop-adjunction2}. Then the linear system $|5K_Y+ 3\det
\mathcal{V}|$ contains an effective nontrivial divisor.
\end{cor}
\begin{proof} Due to our assumptions we know that $A$ is very
ample and $h^0(K_Y+A)\geq 1$. Let $C \in |A|$ be a smooth curve.
Tensoring the exact sequence
$$ 0 \to -A \to \mathcal{O}_Y \to \mathcal{O}_C \to 0 $$
by $2(K_Y + A) +A$ we get the exact sequence
$$0 \to 2K_Y+2A \to 5K_Y+ 3 \det \mathcal{V} \to (5K_Y+ 3 \det \mathcal{V})_C
\to 0.$$ Therefore
$$h^0(5K_Y+ 3\det \mathcal{V}) \geq h^0\big(2(K_Y + A)\big) \geq h^0(K_Y +
A)\geq 1.$$ Suppose that $5K_Y+ 3\det \mathcal{V}$ is trivial.
Then $Y$ is a Del Pezzo surface and $-K_Y$ is divisible by $3$ in
the Picard group. This means that $Y=\mathbb{P}^2$ and $\det
\mathcal{V} = \frac{5}{3}(-K_Y) = \mathcal{O}_{\mathbb{P}^2}(5)$.
But then $(Y, \mathcal{V})$ would be as in case (2) of Proposition
\ref{prop-adjunction2}, a contradiction. Therefore $|5K_Y+ 3\det
\mathcal{V}|$ contains an effective nontrivial divisor.
\end{proof}

\par
Let us analyze more closely the various exceptions arisen in the
discussion. In cases (1)--(4) of Proposition
\ref{prop-adjunction}, clearly $5K_Y+ 3\det \mathcal{V}$ is not
effective. The same is true for case (1) of Proposition
\ref{prop-adjunction2}, while, as already observed, $5K_Y+ 3\det
\mathcal{V}$ is trivial in case (2) of Proposition
\ref{prop-adjunction2}. To complete the discussion it remains to
consider case (5) of Proposition \ref{prop-adjunction} and case
(3) of Proposition \ref{prop-adjunction2}. In both cases, for
every fiber $f$ of the projection $p:Y \to B$ we know that $K_Y\
f=-2$ and $\deg \mathcal{V}_f = 2$ or $3$, respectively. Hence
$(5K_Y + 3 \det \mathcal{V})\ f < 0$ in both cases. Therefore
$5K_Y + 3 \det \mathcal{V}$ is not effective. Let us add something
more on these two cases. First of all we can write
$Y=\mathbb{P}(\mathcal{F})$ for an ample vector bundle
$\mathcal{F}$ of rank $2$ on $B$. Then $K_Y= -2\xi+p^*(K_B+ \det
\mathcal{F})$ where $\xi$ is the tautological line bundle on $Y$.

\par In case (5) of Proposition \ref{prop-adjunction},
since $\mathcal{V} \otimes \xi^{-1}$ restricts trivially to every
fiber $f$ of $p$, there exists a vector bundle $\mathcal{G}$ of
rank $2$ on $B$ such that $\mathcal{V}=\xi \otimes
p^*\mathcal{G}$. Thus $\det \mathcal{V}= 2\xi + p^* \det
\mathcal{G}$ and $c_2(\mathcal{V})= \xi^2 + \xi p^* \det
\mathcal{G}= \deg \mathcal{F} + \deg \mathcal{G}$ by the Chern--Wu
relation. Thus
$$0 < c_1(\mathcal{V})^2 = 4 (\deg \mathcal{F} + \deg \mathcal{G}) = 4
c_2(\mathcal{V}).$$

\par
In case (3) of Proposition \ref{prop-adjunction2} we have
$(\mathcal{V} \otimes [-2\xi])_f = \mathcal{O}_{\mathbb{P}^1}
\oplus \mathcal{O}_{\mathbb{P}^1}(-1)$ for every fiber $f$ of $p$.
Hence $p_*(\mathcal{V} \otimes [-2\xi])$ is a line bundle, say
$\mathcal{A} \in \text{Pic}(B)$, and we have an obvious injection
giving rise to an exact sequence
$$0 \to 2\xi+p^*\mathcal{A} \to \mathcal{V} \to \mathcal{Q} \to
0,$$ where $\mathcal{Q}$ is a very ample line bundle on $Y$, since it is a
quotient of $\mathcal{V}$. This gives $\det \mathcal{V}=2\xi+p^*
\mathcal{A} + \mathcal{Q}$, and from the relation $3 = \deg
\mathcal{V}_f =(2\xi+p^* \mathcal{A} + \mathcal{Q}) f = 2 +
\mathcal{Q} f$, we conclude that $\mathcal{Q}= \xi + p^*
\mathcal{M}$ for some line bundle $\mathcal{M} \in \text{Pic}(B)$.
Thus $\det \mathcal{V}=3\xi+p^*(\mathcal{A}+\mathcal{M})$. Note
also that
$$c_2(\mathcal{V})= 2\xi^2 + \xi\ p^*(\mathcal{A}+2\mathcal{M})=
2\deg \mathcal{F} + \deg \mathcal{A} + 2\deg \mathcal{M}$$ by the
Chern--Wu relation. Thus
$$c_1(\mathcal{V})^2 = 9\deg \mathcal{F} + 6(\deg \mathcal{A}+\deg
\mathcal{M}) = c_2(\mathcal{V}) + (7\deg \mathcal{F} + 5\deg
\mathcal{A}+ 4\deg \mathcal{M}).$$

\medskip
\par The previous discussion allows us to provide a significant
lower bound for $\deg \Phi_2$ of threefold scrolls in
$\mathbb{P}^8$, apart from very restricted cases.

\begin{thm}\label{thm-details}
Let $X \subset \mathbb{P}^{8}$ be a threefold scroll of degree $d$
over a smooth surface $Y$, with projection $\pi$, and let
$\mathcal{V}=\pi_*\mathcal{L}$ where $\mathcal{L}$ is the
hyperplane line bundle. Suppose that $X$ satisfies our general
assumptions. Then $\deg \Phi_2 \geq 3d$ unless
\begin{enumerate}
    \item $(Y, \mathcal{V})=(\mathbb{P}^2,
    \mathcal{O}_{\mathbb{P}^2}(2) \oplus
    \mathcal{O}_{\mathbb{P}^2}(1))$, where $(d,\deg \Phi_2)=(7,3)$;
    \item $(Y, \det \mathcal{V})= (\mathbb{P}^2,
    \mathcal{O}_{\mathbb{P}^2}(4))$, in which case $\deg \Phi_2 = 3d -12$;
    \item $Y=\mathbb{P}(\mathcal{F})$ for some ample
    vector bundle $\mathcal{F}$ of rank $2$ over a smooth curve $B$ of genus
    $q$, and $\mathcal{V}=\xi \otimes p^*\mathcal{G}$, where $\xi$
    is the tautological line bundle, $p\colon Y \to B$ is the bundle projection,
    and $\mathcal{G}$ is a vector bundle of rank $2$ on $B$; in
    this case
    $$\deg \Phi_2 = 3d + 20(q-1) + 2(\deg \mathcal{F} + \deg \mathcal{G});$$
    \item $Y=\mathbb{P}(\mathcal{F})$ for some ample
    vector bundle $\mathcal{F}$ of rank $2$ over a smooth curve $B$ of genus
    $q$, and $\mathcal{V}$ fits into an exact sequence
    $$0 \to 2\xi+p^*\mathcal{A} \to \mathcal{V} \to \xi +p^*\mathcal{M} \to
    0, $$
    where $\xi$ is the tautological line bundle, $p\colon Y \to B$ is the bundle
    projection, and $\mathcal{A}$ and $\mathcal{M}$ are line bundles on $B$;
    in this case $$\deg \Phi_2 = 3d + 30(q-1) + 12\deg \mathcal{F}
    + 8(\deg \mathcal{A}+\deg \mathcal{M}).$$
\end{enumerate}
Moreover, equality $\deg \Phi_2 = 3d$ holds if and only if $(Y,
\det \mathcal{V}) = (\mathbb{P}^2,
\mathcal{O}_{\mathbb{P}^2}(5))$.
\end{thm}
\begin{proof} Dotting formula (\ref{divisor}) with $L^2$ and
taking into account Corollary \ref{cor-effectivity}, we see that
$\deg \Phi_2 = L^2 \Phi_2 = 3L^3 +
c_1(\mathcal{V})(5K_Y+3c_1(\mathcal{V}))
> 3d$, provided that $(Y, \mathcal{V})$ is not any of the exceptions
listed in Propositions \ref{prop-adjunction} and
\ref{prop-adjunction2}, while $\deg \Phi_2 = 3d$ in case (2) of
Proposition \ref{prop-adjunction2}. Cases (5)  of Proposition
\ref{prop-adjunction} and (3) of Proposition
\ref{prop-adjunction2} lead to (3) and (4) in the statement
respectively, and the precise expression for $\deg \Phi_2$ simply
follows from the previous discussion of these cases. Exception (2)
in Proposition \ref{prop-adjunction} gives (1) in the statement;
note that this is the pair discussed in Example 5 for $m=2$.
Exception (1) in Proposition \ref{prop-adjunction2} gives (2).
Finally, exceptions (1), (3) and (4) in Proposition
\ref{prop-adjunction} do not satisfy the assumption on the generic
rank of the second jet map, because $h^0(\mathcal
L)=h^0(\mathcal{V}) \leq 8$ in all these cases.
\end{proof}

\medskip
A similar result can be obtained for fourfold scrolls over
surfaces. In this case the expression provided by Theorem
\ref{thm-divisor} becomes
\begin{equation*}
[\Phi_2] = \pi^*\big(6K_Y + 3c_1(\mathcal{V})\big) + 3L,
\end{equation*}
so that we have to look at the effectivity of $2K_Y+ \det
\mathcal{V}$, where now $\mathcal{V}$ is of rank $3$. Letting $A:=
K_Y+ \det \mathcal{V}$ as before, the analogue of Proposition
\ref{prop-adjunction} is simply that $A$ is very ample unless
$(Y,\mathcal{V}) = (\mathbb{P}^2,
\mathcal{O}_{\mathbb{P}^2}(1)^{\oplus 3})$. Apart from this case,
arguing as in Proposition \ref{prop-adjunction2}, we obtain that
$2K_Y+ \det \mathcal{V}$ is effective unless $(Y, \det
\mathcal{V})$ is as in cases (1), (2) of Proposition
\ref{prop-adjunction2}, or $Y$ is a $\mathbb{P}^1$-bundle over a
smooth curve $B$ and $\mathcal{V}_f =
\mathcal{O}_{\mathbb{P}^1}(1)^{\oplus 3}$ for every fiber $f$ of
the bundle projection $Y \to B$. This immediately leads to a
result of the same type as Theorem \ref{thm-details}.

\medskip
Consider again a threefold scroll $X=\mathbb{P}(\mathcal{V})
\subset \mathbb{P}^8$ satisfying our general assumptions, and let
$g$ be its sectional genus. Then (\ref{two-expressions}) gives
\begin{equation} \label{deg in P8}
\deg \Phi_2 = 10(g-1) + c_1(\mathcal{V})^2 - 3c_2(\mathcal{V}).
\end{equation} Note that for the scroll $X =
\mathbb{P}(\mathcal{O}_{\mathbb{P}^2}(1) \oplus
\mathcal{O}_{\mathbb{P}^2}(2))$ discussed in Example 5 for $n=3$,
the value of $c_1(\mathcal{V})^2 - 3c_2(\mathcal{V})$ is $3$.
Putting some restriction on the vector bundle $\mathcal{V}$, we
get the following characterization.

\begin{prp}
Let $X \subset \mathbb{P}^{8}$ be a threefold scroll over a
smooth surface $Y$, satisfying our general assumptions and such
that $c_1(\mathcal{V})^2 \geq 3c_2(\mathcal{V})$. Then either $X =
\mathbb{P}(\mathcal{O}_{\mathbb{P}^2}(1) \oplus
\mathcal{O}_{\mathbb{P}^2}(2))$ or $\deg \Phi_2 \geq 10(g-1) \geq
20$.
\end{prp}
\begin{proof} (\ref{deg in P8}) gives $\deg \Phi_2 \geq
10(g-1)$, due to the assumption on the Chern classes of
$\mathcal{V}$. Now we use classification results for projective
manifolds with low sectional genus \cite{Io}. Case $g=0$ cannot
occur. Actually, by a result of Sato \cite[Theorem A]{Sa} the only
possibility would be that $X$ is the Segre product
$\mathbb{P}^2\times \mathbb{P}^1 \subset \mathbb{P}^{5}$, which,
as noted in Example 2, does not satisfy our general assumptions.
Suppose $g=1$. Then either $(X,\mathcal{L})$ is a Del Pezzo
threefold, or a scroll over a smooth curve $W$ of genus $1$. The
former case leads to the exception in the statement. Recall that
this $X$ is isomorphic to $\mathbb{P}^3$ blown-up at one point and
$\mathcal{L}=\sigma^* \mathcal{O}_{\mathbb{P}^2}(2) \otimes
\mathcal{O}_X(-E)$, where $\sigma: X \to \mathbb{P}^3$ is the
blowing-up and $E$ is the exceptional divisor. The latter case
cannot occur: actually the fibers of the the scroll over $W$ map
surjectively to the smooth surface $Y$ via the scroll projection
$\pi$, so that $Y=\mathbb{P}^2$ \cite[Example 1.4, p.\ 136]{BPV};
hence $h^1(\mathcal{O}_X)=0$, while, due to the scroll structure
over $W$, it should be $h^1(\mathcal{O}_X)=1$. Let $g=2$. As in the
previous case we exclude that $(X,\mathcal{L})$ is a scroll over a
smooth curve of genus $2$. Then \cite[Corollary 3.3]{Io}
$(X,\mathcal{L})$ is a quadric fibration over $\mathbb{P}^1$ and,
using \cite[Theorem, p.~4]{L}, the scroll structure over $Y$ allows us to
recognize $X$ as the general projection in $\mathbb{P}^8$ of the
Segre product $\mathbb{F}_1\times \mathbb{P}^1 \subset
\mathbb{P}^9$, where $\mathbb{F}_1 \subset \mathbb{P}^4$ is the
rational cubic scroll \cite[Theorem 3.4, case vi)]{Io}. However,
as observed in Example 2, this $X$ does not satisfy the assumption
$s_2=9$. Therefore $g \geq 3$ and the assertion is proved.
\end{proof}

Let $Y \subset \mathbb{P}^5$ by any generically 2-regular surface.
Using (\ref{divisor}) and Example 1, we can provide an application
to the general projection $X \subset \mathbb{P}^8$ of the product
scroll $X^o:=Y \times \mathbb{P}^1 \subset \mathbb{P}^{11}$.
According to \cite[Proposition 0.3]{Sh}, the class of $\Phi_2(Y)$
is $4K_Y+6H$, $H$ being the hyperplane class. Hence $\Phi_2(X^o) =
\pi^{-1}(\Phi_2(Y))$ according to Example 1. This means that
$\Phi_2(X^o)$ is a divisor. Projecting generically $X^o$ to
$\mathbb{P}^8$, we expect that the image $X \cong X^o$ has the
whole $\mathbb{P}^8$ as osculating space at the general point.
Thus $\Phi_2:=\Phi_2(X)$ is a divisor containing the projection
$\Phi_2^o$ of $\Phi_2(X^o)$ as a component, plus some other
components, i.\,e., $\Phi_2 = \Phi_2^o + R$, where $R$ is an
effective divisor. In fact, our formula (\ref{divisor}) says that
$[\Phi_2] = 3L + \pi^*(5K_Y+6H)$, since $\mathcal{V}=H^{\oplus
2}$. As $\Phi_2^o = \pi^*(4K_Y+6H)$, this means that
$R=3L+\pi^*K_Y$. Note that the effectiveness of $R$ is independent of
$(Y,H)$. Indeed, $L= Y_0+\pi^*H$, where $Y_0$ is the
tautological section of the trivial bundle $\mathcal{O}_Y^{\oplus
2}$. Hence $R= 3Y_0 + \pi^*(K_Y+3H)$ and $K_Y+3H$ is effective,
$Y$ being generically $2$-regular. Since $X$ was obtained by
projecting $X^o$ from a general $\mathbb{P}^2 \subset
\mathbb{P}^{11}$, the above conclusion can be converted into an
information on $X^o$. In particular, it says that the part of
$\Phi_2$ coming from the osculating spaces to $X^o \subset
\mathbb{P}^{11}$ which meet a general plane, is a surface of degree
$(3L+\pi^*K_Y) L^2 = 3d + 2HK_Y$.

\section{The case $k=2$: formulas for threefold scrolls over surfaces}
\label{k=2}
Let $X \subset \mathbb{P}^N$ be a $n$-dimensional
scroll over a smooth surface $Y$, satisfying our general
assumptions. Then (\ref{explicit}) gives (\ref{tag
8}), that we rewrite for the convenience of the reader:
\begin{equation}
 c(\mathcal E_2)= \pi^* c(\mathcal V^{\vee})\ \pi^* c(\mathcal V^{\vee}
\otimes T_Y)\ c(\pi^*S^2 T_Y \otimes \mathcal L^{-1}). \nonumber
\end{equation}
According to Theorem \ref{thm-2}, recall that to determine the
class of $\Phi_2$ we have to take the inverse of the expression
in (\ref{tag 8}). Recall that here both $\mathcal V^{\vee}$ and
$\mathcal V^{\vee} \otimes T_Y$ are vector bundles over a surface.
Hence,
\begin{equation} \label{tag 9}
c(\mathcal V^{\vee})^{-1} = 1+ v_1 + [v_1^2-v_2],
\end{equation}
where, for shortness, we have set $v_i = c_i(\mathcal{V})$. On the
other hand, let $c_i = c_i(T_Y)$. Then computing the Chern classes
of $\mathcal V^{\vee} \otimes T_Y$ we get
\setlength\arraycolsep{2pt}
\begin{eqnarray*} \textstyle c_1(\mathcal V^{\vee} \otimes T_Y)&=& -2v_1+(n-1)c_1\\\textstyle
c_2(\mathcal V^{\vee} \otimes T_Y)& = &v_1^2+2v_2 - (2n-3) v_1 c_1
 +\textstyle  \binom{n-1}2 c_1^2 + (n-1)c_2.
\end{eqnarray*}
Thus,
\setlength\arraycolsep{2pt}
\begin{eqnarray} \label{tag 10}\textstyle
c(\mathcal V^{\vee} \otimes T_Y)^{-1}& = &1 + \big[2v_1
- (n-1)c_1\big]  + \big[3v_1^2
 - 2v_2 \nonumber \\ &&- (2n-1) v_1 c_1
+ \textstyle  \binom{n}2 c_1^2 - (n-1)c_2\big].
\end{eqnarray}
As to the last factor in (\ref{tag 8}), note that $\pi^*S^2T_Y
\otimes \mathcal L^{-1}$ is a vector bundle of rank $3$ on $X$.
First of all, the vector bundle $S^2T_Y$ has total Chern class
$$c(S^2T_Y) = 1 + 3c_1 + \big(2c_1^2 +
4c_2\big).$$ So, letting $L=c_1(\mathcal L)$ and
$C_i=\pi^*c_i=\pi^*c_i(T_Y)$, we obtain
\setlength\arraycolsep{2pt}
\begin{eqnarray} \label{tag 11}
c(\pi^*S^2T_Y \otimes \mathcal L^{-1}) & =& 1  + \big[3
C_1 -3L \big]\nonumber\\
&& +\big[2C_1^2 + 4C_2 - 6 C_1L + 3L^2 \big]
\nonumber\\
&& -\big[(2C_1^2 + 4C_2)L - 3C_1 L^2 + L^3\big].
\end{eqnarray}
Computing the inverse of this total Chern class is very difficult
in general. Let us do it for threefold scrolls. So from now on in
this section we let $n=3$. Recall that for a vector bundle $B$ of
rank $3$ with $b_i= c_i(B)$ we have
$$c(B)^{-1}=1 - b_1 + [b_1^2-b_2] +[-b_1^3 + 2b_1b_2 -
b_3].$$
\medskip
Thus, if $n=3$, the inverse of (\ref{tag 11}) is:
\setlength\arraycolsep{2pt}
\begin{eqnarray*}
c(\pi^*S^2T_Y \otimes \mathcal L^{-1})^{-1}&=& 1 + 3
\big[L - C_1\big] \nonumber\\ && + \big[6L^2 + 7C_1^2 -
12C_1L-4C_2\big]\nonumber\\ && + 5\big[2L^3 + 7C_1^2L - 4C_2L
-6C_1L^2\big].
\end{eqnarray*}

Letting $V_i = \pi^*v_i= \pi^*c_i(\mathcal V)$, the product of the
three factors in (\ref{tag 8}) gives the following expression for
\setlength\arraycolsep{2pt}
\begin{eqnarray}\label{tag 12}
c(\mathcal{E}_2)^{-1} &=& 1 + \big[3L + 3V_1 -5C_1\big] +
\big[6L^2 + 9 V_1 L- 18C_1L \nonumber\\
&&+ 6 V_1^2 - 3V_2 -16
C_1V_1 + 16C_1^2 -6C_2\big] \nonumber\\
&&+\big[10L^3
-42C_1L^2 + 18 V_1L^2 +68 C_1^2L \nonumber\\
&&-26 C_2L -57 C_1V_1L + 18 V_1^2L - 9V_2L \big].
\end{eqnarray}
In particular, if $\big[\Phi_2\big]$ has codimension $l=n=3$, then
$X$ has finitely many flexes, their number being given by
\begin{equation} \label{tag 13}
\deg\Phi_2=19d+68c_1^2-26c_2-99c_1v_1 +27v_1^2 ,
\end{equation}
which fits with the computations already done when $Y$ is an
abelian surface (\ref{abelian}).
\medskip

We shall now apply the formulas to  threefold scrolls $X$ over some
special surfaces $Y$. According to (\ref{range,k=2}) we should consider $X
\subset \mathbb P^N$ with $N=10,9,8$, but in view of the
discussion in Section \ref{k=m=2,divisor}, we need only consider
the cases $N=10$ and $9$. For some of the needed computations  it is
better to recall that $X=\mathbb{P}(\mathcal{V})$, where
$\mathcal{V}=\pi_*\mathcal{L}$ is a very ample vector bundle of
rank $2$ on $Y$. Then $K_X=-2\mathcal L+\pi^*(K_Y+\det \mathcal V)$, where
$\pi:X \to Y$ is the scroll projection. Moreover, $L^2= L\, V_1 -
V_2$, where $V_i = \pi^*c_i(\mathcal{V})$, by the Chern--Wu
relation. In particular, this gives $L^2\,\pi^*D = L \,V_1\, \pi^*D =
c_1(\mathcal{V})\,D$ for any divisor $D$ on $Y$.

\medskip
First of all consider $X:=
\mathbb{P}(\mathcal{O}_{\mathbb{P}^2}(2)^{\oplus 2}) \subset
\mathbb{P}^{11}$. Note that $X$ is the Segre product of the
Veronese surface in $\mathbb{P}^5$ with $\mathbb{P}^1$. Hence
$\Phi_2(X) = \emptyset$ by what we said at the end of Example 1
combined with \cite[Proposition 3.8]{Sh}. This could lead to doubts
concerning formula (\ref{tag 13}). However, note that here $N=11$
instead of $10$. So the formula does not apply to $X$, but it
applies to its general projection in $\mathbb{P}^{10}$, for which
we get $\deg \Phi_2 = 6$. Rephrasing this in terms of the linearly
normal scroll $X \subset \mathbb{P}^{11}$ gives:

\begin{cor}
Through a general point of $\mathbb{P}^{11}$ there pass $6$ second order
osculating spaces of the Segre product $S \times \mathbb{P}^1$,
where $S$ is the Veronese surface in $\mathbb{P}^5$.
\end{cor}

\medskip
A nice result concerning scrolls $X \subset \mathbb{P}^{10}$,
derived from formula (\ref{tag 13}), is the following:

\begin{prp} \label{uninflectedinP10}
Let $X \subset \mathbb{P}^{10}$ be a threefold scroll over
$\mathbb{P}^2$, satisfying our general assumptions. Then
$\Phi_2(X) \not= \emptyset$.
\end{prp}

\begin{proof} If $X$ is uninflected, then Theorem \ref{thm-2} applies, and
hence the expression in formula (\ref{tag 13}) is equal to $0$. As $Y=\mathbb{P}^2$, letting
$(x,y)=\big(c_1(\mathcal{V}),c_2(\mathcal{V})\big)$, (\ref{tag
13}) becomes $\deg \Phi_2 = 46x^2-297x+534-19y$. Then $(x,y)$ must
be be an integral point of the parabola $C \subset \mathbb{A}^2$
defined by $y=\frac{1}{19}(46x^2-297x+534)$. On the other hand,
the point $(x,y)$ must lie in the region $R \subset \mathbb{A}^2$
defined by $y \leq x^2-8$, since $x^2-y = d \geq
\text{codim}_{\mathbb{P}^{10}}(X) + 1=8$. But an immediate check
shows that there are no integral points on the arc $R \cap C$.
\end{proof}

\medskip
In the following we assume that our threefold scrolls $X \subset
\mathbb{P}^N$ satisfy the general assumptions of the discussion
leading to the formulas in Section 4.

\medskip Now consider scrolls $X \subset \mathbb{P}^9$.
Here $\ell=2$, so that
$\Phi_2$ is a $1$-cycle, whose class is the second term in
(\ref{tag 12}). Dotting it with $L$ and taking into account the
Chern--Wu relation, we get the following expression for its
degree:
\begin{equation} \label{deg in P9}
\deg \Phi_2 = 9d + 12c_1(\mathcal{V})^2+34 K_Yc_1(\mathcal{V}) +
16K_Y^2 - 6 c_2(T_Y).
\end{equation}

First suppose $Y = \mathbb{P}^2$. In this case, letting
$v=c_1(\mathcal{V})$, (\ref{deg in P9}) becomes
\begin{equation} \label{deg in P9 forY=P2}
\deg \Phi_2 = 9d + 6v(2v-17)+126.
\end{equation}

This allows us to prove an analogue of Proposition
\ref{uninflectedinP10}.

\begin{prp} \label{uninflectedinP9}
Let $X \subset \mathbb{P}^{9}$ be a threefold scroll over
$\mathbb{P}^2$, satisfying our general assumptions. Then
$\Phi_2(X) \not= \emptyset$.
\end{prp}

\begin{proof} If $\Phi_2(X)= \emptyset$ then $\deg \Phi_2 =0$, so
that (\ref{deg in P9 forY=P2}), after dividing by $3$, gives
$$3(d+14) = 2v(17-2v).$$
Note that $X$ has degree $d \geq \text{codim}_{\mathbb{P}^{9}}(X)
+ 1=7.$ Hence the left hand term is $\geq 63$. Then looking at the
right hand term we conclude that $v$ can only be $3,4,5$ or $6$.
However, $v=5$ is not possible, otherwise the right hand term would
not be divisible by $3$. If $v=3$ then the equation above gives
$d=8$, hence $8= v^2 - c_2(\mathcal{V})= 9 - c_2(\mathcal{V})$,
i.\,e., $c_2(\mathcal{V})=1$, but this would imply that
$\mathcal{V} =\mathcal{O}_{\mathbb{P}^2}(1)^{\oplus 2}$
\cite[Theorem 11.1.3] {BS}, which is clearly not compatible with
$v=3$. If $v=6$, then the equation above gives $d=6$. On the
other hand, recalling what we said in the proof of Proposition
\ref{Theorem-abelian}, the curve genus of $(X,\mathcal{L})$ is
$g=g(\mathbb{P}^2,\mathcal{O}_{\mathbb{P}^2}(v))=10$, since $v=6$.
But this implies that the general curve section of $X$ is a plane
curve, which is clearly impossible. The only remaining possibility
is $v=4$, which gives $g=3$, $d=10$ and $c_2(\mathcal{V})=6$.
These invariants characterize the Bordiga scrolls of degree $10$,
but, as observed in Example 3, they do not satisfy our general
assumptions.
\end{proof}

\medskip
The next result takes care of the case $Y=\mathbb{F}_e$, a
Segre--Hirzebruch surface.

\begin{prp} \label{F_e}
Let $X \subset \mathbb{P}^{9}$ be a threefold scroll over
$\mathbb{F}_e$ $(e \geq 0)$. Then $\Phi_2(X) \not= \emptyset$,
except possibly if $\mathcal{V}$ is uniform of type $(1,1)$ on the
fibers $\mathbb{F}_e$, and
$$9d - 32 = 20(b-e),$$ where $b$ is the degree of $\mathcal{V}$
restricted to a minimal section of $\mathbb{F}_e$.
\end{prp}
\begin{proof} Let $s$ and $f$ be a section of minimal
self-intersection ($s^2 =-e$) and a fiber of $\mathbb{F}_e$
respectively. Since the classes of $s$ and $f$ generate
$\text{Pic}(\mathbb{F}_e)$, we can write $\det \mathcal V=as +
bf$, where $a = \deg \mathcal{V}_f \geq \text{rk}\mathcal{V} = 2$,
equality implying that $\mathcal{V}_f =\mathcal{O}_{\mathbb
P^1}(1)^{\oplus 2}$ for every fiber $f$. Moreover $b-ae =
(as+bf)s =\deg \mathcal{V}_s \geq 2$. Recalling that
$K_{\mathbb{F}_e}=-2s -(2+e)f$, (\ref{deg in P9}) gives $\deg
\Phi_2 = 9d + 12(2b-ae)a + 34(ae - 2a-2b) + 104$. So, if
$\Phi_2(X)= \emptyset$, then $(a,b)$ must be an integral point in
the region $R$ of the affine plane defined by $a \geq 2, b \geq
ea+2$, lying on the hyperbola $\gamma$ with equation
\begin{equation} \label{hyperbola}
12ea^2 - 24ab  + 34(2-e)a +68b - (9d+104)=0.
\end{equation}
Let us draw a picture. The asymptotes of $\gamma$ are the lines
parallel to the $a$-axis and to the line of equation $b =
\frac{e}{2}a$ passing through the center $A$, which has
coordinates $(\frac{17}{6},\frac{17}{6}+\frac{17}{12}e)$. Note
that the second asymptote cuts the line $a=2$ at the point $T =
(2, e+\frac{17}{6})$. On the other hand, the hyperbola $\gamma$
cuts the line $a=2$ at the point $B = (2, e + \frac{9d-32}{20})$.
Note that as a consequence of (\ref{hyperbola}), $d$ must be even.
On the other hand, $d \geq \text{codim}_{\mathbb{P}^9}(X)+1=7$, as
already observed. Hence either $d=8$ or $d \geq 10$. Suppose that
$d=8$. Then $(X,\mathcal{L})$ has $\Delta$-genus $\leq 1$ (since
our $X \subset \mathbb{P}^9$ is not necessarily linearly normal).
This implies that it is either a scroll over $\mathbb{P}^1$ or a
Del Pezzo threefold \cite[Proposition 2.4]{Io}. The former case
cannot occur because our $X$ has Picard number $3$. The latter
case also gives a contradiction, since the only Del Pezzo
threefold of degree $8$ is the pair $(\mathbb{P}^3,
\mathcal{O}_{\mathbb{P}^3}(2))$. Therefore $d \geq 10$. This fact
implies that the second coordinate of $B$ is bigger than that of
$T$. Now, looking at the position of the region $R$ we can see
that $\gamma \cap R$ does not contain integral points on the
branch lying in the half plane $a > 2$. As a consequence, the
integral points on $\gamma \cap R$ we are looking for are confined
between the line $a=2$ and the vertical asymptote
$a=\frac{17}{6}$. Since this value is smaller than $3$ it turns
out that they consist of the single point $B$. Thus $a=2$, which
implies that $\mathcal{V}$ is uniform of type $(1,1)$ on the
fibers $f$, and $b=e + \frac{9d-32}{20}$, which leads to the
numerical condition in the statement.
\end{proof}

\medskip
By the same method we can produce a result taking care of another
interesting class of surfaces, namely products of curves. Since
the case $Y = \mathbb{F}_0 = \mathbb{P}^1 \times \mathbb{P}^1$ has
already been considered, we will assume that not both factors are
rational curves.

\begin{prp} \label{products}
Let $X \subset \mathbb{P}^{9}$ be a threefold scroll over a surface
$Y$ which is the product of two smooth curves, not both rational. Then
$\Phi_2(X) \not= \emptyset$, except possibly for the following case:
$Y=B \times \mathbb P^1$, where $B$ has genus $q \geq 1$,
$\mathcal{V}$ is uniform of type $(1,1)$ on the fibers of the
first projection, and
$$9d + 32 (q-1) = 20 b,$$ where $b$ is the degree of $\mathcal{V}$
restricted to the fibers of the second projection (in particular
$b \geq 5$).
\end{prp}
\begin{proof} Set $Y = C_1 \times C_2$ and let $g_i$ be the genus
of the curve $C_i$, $i=1,2$. Up to renaming we can suppose that
$g_1 \geq g_2$. Since $K_Y^2 = 8(g_1-1)(g_2-1)$,
$c_2(T_Y)=4(1-g_1)(1-g_2)$, and $\deg \mathcal{V}_{C_i}>0$, (\ref{deg
in P9}) immediately implies that $\deg \Phi_2 \geq 9d + 12
c_1(\mathcal{V})^2 > 0$ if $g_2 \geq 1$. So, we can put $C_2
=\mathbb P^1$. In other words, it only remains to consider the
case when $Y= B \times \mathbb P^1$, where $B$ has genus $q \geq
1$. Letting $f$ and $s$ denote  the fibers of the first and the second
projection respectively, $\det \mathcal{V}$ is numerically
equivalent to $as+bf$, where $a = \deg \mathcal{V}_f \geq
\text{rk}\,\mathcal{V} = 2$, equality implying that $\mathcal{V}_f
=\mathcal{O}_{\mathbb P^1}(1)^{\oplus 2}$ for every fiber $f$.
Moreover $b =\deg \mathcal{V}_s \geq 5$, since there are no
irrational surface scrolls of degree $<5$. Then (\ref{deg in P9})
gives $\deg \Phi_2 = 9d + 24 ab + 68(q-1)a - 68b + 104(1-q)$. So,
if $\Phi_2(X)= \emptyset$, then $(a,b)$ must be an integral point
in the region $R$ of the affine plane defined by $a \geq 2, b\geq
5$, lying on the hyperbola $\gamma$ with equation
$$24ab + 68(q-1)a -68b + 9d-104(q-1)=0.$$
The asymptotes of $\gamma$ are the parallel lines to the $a$-axis
and to the $b$-axis passing through the center $A$, which has
coordinates $(\frac{17}{6},\frac{17}{6}(1-q))$. On the other hand,
$\gamma$ cuts the line $a=2$ at the point $B$ whose second
coordinate is $\frac{1}{20}(9d + 32(q-1))$, which is bigger than
that of $A$. Then the integral points on $\gamma \cap R$ we are
looking for are confined between the line $a=2$ and the vertical
asymptote, and since the abscissa of $A$ is smaller than $3$, they
consist of the single point $B$. Thus $a=2$, which implies that
$\mathcal{V}$ is uniform of type $(1,1)$ on the fibers $f$, and
$b=\frac{1}{20}(9d + 32(q-1))$, which leads to the numerical
condition in the statement.
\end{proof}

\medskip
If $Y$ is an abelian surface, then (\ref{deg in P9}) gives $\deg
\Phi_2 = 9d + 12(2g-2)$, in accordance with the general formula in
Proposition \ref{Theorem-abelian}. For $Y$ a K3 surface, (\ref{deg
in P9}) gives $\deg \Phi_2 = 3(3d + 4(2g-2)-48)$. Recall that $d
\geq 7$ as already observed. Let $S$ be a general surface section
of $X$ and look at the reduction $(Y, \det \mathcal{V})$ of
$(S,\mathcal{L}_S)$. The general curve $C \in |\det \mathcal{V}|$
is a canonical curve of degree $d + c_2(\mathcal{V}) \geq 8$,
hence $g(C) \geq 5$. Therefore $g \geq 5$ and then $\deg \Phi_2
\geq 15$. Though this is a rough inequality, it is enough to
conclude that:

\begin{prp}
Let $X \subset \mathbb{P}^{9}$ be a threefold scroll over a K$3$
surface. Then $\Phi_2(X) \not= \emptyset$.
\end{prp}

We can prove a similar result for $Y$ a minimal surface of Kodaira
dimension $1$, in two cases: when the base curve of the elliptic
fibration has positive genus, and when the base is $\mathbb P^1$
provided that $Y$ is an elliptic quasi-bundle in the sense of
Serrano \cite[Definition 1.1 and Lemma 1.5]{Se}.

\section{The case $k=2$: formulas for fourfold scrolls over threefolds}
\label{k=2,m=3}

Let $X \subset \mathbb{P}^N$ be a fourfold scroll over a smooth
threefold $Y$, satisfying our general assumptions.
Again, according to Theorem \ref{thm-2}, to determine
the class of $\Phi_2$ we have to take the inverse of the
expression in (\ref{tag 8}). Now both $\mathcal V^{\vee}$ and
$\mathcal V^{\vee} \otimes T_Y$ are vector bundles over a
threefold. Hence (\ref{tag 9}) has to be replaced with
\begin{equation*}
c(\mathcal V^{\vee})^{-1} = 1 + v_1 + [v_1^2 - v_2] +
[v_1^3 - 2 v_1v_2],
\end{equation*}
where $v_i = c_i(\mathcal{V})$. In the following, as in Section 6
we set also $c_i=c_i(T_Y)$. On the other hand, by using the
splitting principle confined to $\mathcal{V}^{\vee}$, we easily get
the following expressions for the Chern classes of $\mathcal
V^{\vee} \otimes T_Y$. They are:
\setlength\arraycolsep{2pt}
\begin{eqnarray*}&&c_1(\mathcal V^{\vee} \otimes T_Y)= -3v_1+2c_1\\
&&c_2(\mathcal V^{\vee} \otimes T_Y)= 3v_1^2 + 3v_2 -5c_1v_1 + c_1^2
+ 2c_2\\
&&c_3(\mathcal V^{\vee} \otimes T_Y)= -v_1^3 -6 v_1v_2
+4c_1(v_1^2+v_2) - \big( 2c_1^2+4c_2\big) v_1 + 2c_1c_2 + 2c_3.
\end{eqnarray*}
Thus, the analogue of (\ref{tag 10}) for $n=4$ is now the
following:
\setlength\arraycolsep{2pt}
\begin{eqnarray*}
c(\mathcal V^{\vee} \otimes T_Y)^{-1}& = &1 + \big[3v_1
-2c_1\big]\\
& &+ \big[6v_1^2 - 3v_2 - 7c_1v_1 + 3c_1^2
-2c_2\big]\\
&& + \big[10 v_1^3 - 12 v_1v_2 - 16c_1v_1^2 +
8c_1v_2\\
& & +12c_1^2v_1 - 8c_2v_1 - 4 c_1^3 + 6c_1c_2 - 2c_3\big].
\end{eqnarray*}
As to the last factor in (\ref{tag 8}), note that $\pi^*S^2T_Y
\otimes \mathcal L^{-1}$ is a vector bundle of rank $6$ on $X$.
First of all, the vector bundle $S^2T_Y$ has total Chern class
\setlength\arraycolsep{2pt}
\begin{eqnarray}
c(S^2T_Y) = 1 + 4c_1 + 5\big(c_1^2 + c_2\big) +
2c_1^3+11c_1c_2+7c_3. \nonumber
\end{eqnarray}
So, letting $L=c_1(\mathcal L)$ and $C_i=\pi^*c_i=\pi^*c_i(T_Y)$
as in Section 6, we obtain
\setlength\arraycolsep{2pt}
\begin{eqnarray} \label{tag 11'}
c(\pi^*S^2T_Y \otimes \mathcal L^{-1}) & =& 1  + \big[4
C_1-6L \big] \nonumber\\
&&+\big[5\big(C_1^2 + C_2\big)- 20 C_1L
+ 15L^2 \big] \nonumber\\
&& +\big[2C_1^3 + 11C_1C_2+7C_3\nonumber\\
&& -20\big(C_1^2 + C_2\big)L
+ 40C_1L^2 - 20L^3 \big] \nonumber\\
&& +\big[-3\big(2C_1^3+11C_1C_2+7C_3\big)L
\nonumber\\
& &+ 30 \big(C_1^2 + C_2\big)L^2 -40C_1L^3 + 15L^4\big].
\end{eqnarray}

Recall that for a vector bundle $B$ on a fourfold, with $b_i=
c_i(B)$, we have
\setlength\arraycolsep{2pt}
\begin{eqnarray}
c(B)^{-1} & =& 1 - b_1 + [b_1^2-b_2] +[-b_1^3 + 2b_1b_2
- b_3] \nonumber\\
&& +[b_1^4 -3b_1^2b_2 +2b_1b_3 +b_2^2-b_4]. \nonumber
\end{eqnarray}
Thus, since $n=4$, the inverse of (\ref{tag 11'}) is:
\setlength\arraycolsep{2pt}
\begin{eqnarray}
c(\pi^*S^2T_Y \otimes \mathcal L^{-1})^{-1}&=& 1 +
\big[6L - 4C_1\big]\nonumber\\
& &+  \big[21L^2 - 28C_1L +11C_1^2
- 5C_2\big]\nonumber\\
&& +  \big[56L^3 - 112C_1L^2 + 88C_1^2L
-40C_2L \nonumber\\
& & -  26C_1^3 + 29C_1C_2
-7C_3\big] \nonumber\\
&& + \big[126L^4 - 336 C_1L^3 + 396C_1^2L^2 -180C_2L^2
\nonumber\\
& & -  234C_1^3L +261C_1C_2L - 63C_3L \big].\nonumber
\end{eqnarray}
Now, letting $V_i = \pi^*v_i = \pi^*c_i(\mathcal V)$, the product
of the inverses of three factors in (\ref{tag 8}) (which involves
a tedious computation) will give the following expression for
$\Phi_2$, according to whether $\ell=2,3,4$ respectively (for the
case $\ell=1$, see Section \ref{k=2,divisor}).
\setlength\arraycolsep{2pt}
\begin{eqnarray*}
[\Phi_2] & = & 21L^2 + 24V_1L- 40C_1L + 10 V_1^2 - 4V_2
 -25C_1V_1 +  22C_1^2 - 7C_2,
\end{eqnarray*}
\setlength\arraycolsep{2pt}
\begin{eqnarray*}
[\Phi_2]
& = & 56L^3 - 154C_1L^2 + 84V_1L^2 + 162C_1^2L -  52C_2L -166C_1V_1L\\
&& + 60V_1^2L - 24V_2L
 + 20V_1^3 -20V_1V_2 - 65C_1V_1^2 +26C_1V_2\\
&& + 95C_1^2V_1 -30C_2V_1 -64C_1^3 +53C_1C_2 - 9C_3,
\end{eqnarray*}
\setlength\arraycolsep{2pt}
\begin{eqnarray*}
[\Phi_2] & = & 126L^4 + 224V_1L^3 - 448 C_1L^3 - 84V_2L^2
+210V_1^2L^2  -  637C_1V_1L^2\\
&&
+683C_1^2L^2 -222C_2L^2
 +  694C_1^2V_1L -220C_2V_1L
+120V_1^3L\\
&& -  120V_1V_2L - 430C_1V_1^2L
+172C_1V_2L
 -  518C_1^3L\\ && +433C_1C_2L -75C_3L.
\end{eqnarray*}

To compute the degree of $\Phi_2$ in all cases we have to dot the
above expressions with $L^{4-\ell}$. Notice that using the
Chern--Wu relation we can avoid the appearance of the term
$v_1v_2$ in all expressions. Recall that $c_i=c_i(T_Y)$ and
$v_i=c_i(\mathcal{V})$. First of all, if $\Phi_2$ has codimension
$\ell=1$, then $\deg \Phi_2$ is given by (\ref{deg Phi2}), which,
in the present case $(n,m)=(4,3)$ becomes
\begin{equation*}
\deg\Phi_2 = 8d+2v_1^3-6c_1v_1^2 +6c_1v_2.
\end{equation*}
If $\ell=2$, then $\Phi_2$ is a $2$-cycle of degree
\begin{equation*}
\deg\Phi_2 = 35d+20v_1^3-65c_1v_1^2 +40c_1v_2 +22c_1^2v_1-
7c_2v_1.
\end{equation*}
If $\ell=3$, then $\Phi_2$ is a $1$-cycle of degree
\setlength\arraycolsep{2pt}
\begin{eqnarray}
\deg\Phi_2 &=& 120d +100v_1^3 -385c_1v_1^2 +180c_1v_2
+257c_1^2v_1- 82c_2v_1
\nonumber\\
&& - 64c_1^3 +53 c_1c_2 -9c_3.\nonumber
\end{eqnarray}
Finally, if the codimension $\ell$ of $\Phi_2$ is $4$, then $X$
has finitely many flexes, and their  number is given by
\setlength\arraycolsep{2pt}
\begin{eqnarray}
\deg\Phi_2 &=& 340d+340v_1^3-1515c_1v_1^2 + 620c_1v_2 +
1377c_1^2v_1 -442c_2v_1
\nonumber\\
&& -518c_1^3 + 433 c_1c_2 -75c_3.\nonumber
\end{eqnarray}

Note that, when $Y$ is an abelian threefold, the previous
expressions of $[\Phi_2]$ are consistent with the result provided
by Proposition \ref{Theorem-abelian3}.

\medskip

Here is an application of our formulas.

\begin{thm} \label{uninflectedoverP3}
Let $X \subset \mathbb{P}^{N}$ be a fourfold scroll over
$\mathbb{P}^3$, satisfying our general assumptions. Then
$\Phi_2(X) \not= \emptyset$.
\end{thm}

\begin{proof} Clearly, $\Phi_2(X)$ cannot be a divisor, by the Remark
in Section 4. So $\ell \geq 2$. Suppose that $\Phi_2(X)=
\emptyset$; then $\deg \Phi_2 =0$. For $Y=\mathbb{P}^3$ we have
$c_1=c_3=4$, $c_2=6$. Set $x=v_1, y=v_2$. Then, according to the
previous formulas, $\deg \Phi_2$ is given by:
\setlength\arraycolsep{2pt}
\begin{eqnarray*} &&5(7d + 4x^3 -52x^2 + 32y + 62x),\\
&&20(6d + 5x^3 -77x^2 + 36y + 181x -143),\\
&&20(17d + 17x^3 -303x^2 + 124y + 969x -1153),
\end{eqnarray*}
according to whether $\ell = 2,3,4$ respectively. First, consider the
case $\ell=4$. Condition $\deg \Phi_2=0$ gives
\setlength\arraycolsep{2pt}
\begin{eqnarray}\label{a}
17d + 124y = -17x^3 + 303x^2 - 969x + 1153\ .
\end{eqnarray}
Note that the left hand side is positive (since $d > 0$ and
$\mathcal{V}$ is ample), while the right hand side is negative for
$x\gg 0$. In fact, for it to be positive, we need $x \leq 14$. On the
other hand, $x \geq 2$ by Remark A1) in Section 5. Now, by using the
Chern--Wu relation, which in the present case says that $d = x^3 -
2xy$, equation (\ref{a}) can be rewritten as
$$34x^3 - 34xy -303x^2 + 124y + 969x -1153=0 ,$$
which gives
\setlength\arraycolsep{2pt}
\begin{eqnarray}\label{c}
y = \frac{34x^3 -303x^2 + 969x -1153}{2(17x - 62)}\ .
\end{eqnarray}
However, a close check shows that this expression does not produce
a positive integer for any integer $x$ such that $2 \leq x \leq 14$.
Therefore $\ell \not=4$. Case $\ell = 3$ can be ruled out with a
similar procedure. In this case, the analogue of (\ref{a}) is
$$6d + 36y = -5x^3 + 77x^2 - 181x + 143,$$
and the positivity of the left hand side implies that $x \leq 12$.
The analogue of (\ref{c}), noting that $x$ cannot be 3, is
$$y = \frac{11x^3 - 77x^2 + 181x - 143}{12(x - 3)}\ ,$$
and a direct check shows that $y$ is not an integer for $x=2$ and
for any integer $x$ with $4 \leq x \leq 12$. So the only case left is $\ell = 2$.
Here the analogue of (\ref{a}) is
$$7d + 32y = -2x(2x^2-26x+31),$$
and the positivity of the left hand side implies that $x \leq 11$.
On the other hand, the analogue of (\ref{c}) is
$$y = \frac{x(11x^2 - 52x + 62)}{2(7x - 16)}\ .$$
For $x$ in the range $2 \leq x \leq 11$, we can see that this
expression provides a positive integer only for $x=4$ and $8$.
However, for $x=8$, condition $y < \frac{x^2}{2}$, coming from the
Chern--Wu relation, is not satisfied. The case $x=4$ gives
$y=5$. Clearly, if $\mathcal{V}$ is a rank 2 vector bundle with
Chern classes $v_1=4, v_2=5$, then it is indecomposable. Moreover,
since $v_1=4$, the canonical bundle formula shows that $K_X =
-2\mathcal{L}$. Thus $X$ is a Fano fourfold of index $2$. In
particular, $X$ is a Fano bundle in the terminology of \cite{SW},
and since $\mathcal{V}$ is indecomposable, it follows from
\cite[Theorem 2.1]{SW} that $\mathcal V$ is a twist of the null correlation
bundle $N$ fitting into an exact sequence
$$0 \to N \to T_{\mathbb{P}^3}(-1) \to \mathcal{O}_{\mathbb{P}^3}(1)
\to 0$$ \cite[p. 77]{OSS}. In fact, from $(v_1,v_2)=(4,5)$, we
deduce that $\mathcal{V}=N(2)$. Dualizing the above sequence and
using the fact that $N \cong N^{\vee}$, we see that $N(1)$ is a
quotient of $\wedge^2 \big(T_{\mathbb{P}^3}(-1)\big)$, hence
spanned. Therefore, $N(2)=N(1) \otimes
\mathcal{O}_{{\mathbb{P}^3}}(1)$ is very ample. However,
$h^0(N(2))=16$, hence for $X=\mathbb{P}(N(2))$ the tautological
line bundle $\mathcal{L}$ provides an embedding into
$\mathbb{P}^{15}$, so that, if our general assumptions are
satisfied, we should have $\ell=3$ instead of $2$. To have $\ell =
2$, our scroll $X$ would be a general projection in
$\mathbb{P}^{14}$ of this scroll, say $X^o$. However, in this case
$\Phi_2(X)$ could not be empty. Actually, for some $x \in X^o$,
the osculating space to $X^o$ at $x$ contains the center of
projection, hence the image of $x$ in $X$ belongs to $\Phi_2(X)$.
\end{proof}

Now let $Y=\mathbb{Q}^3 \subset \mathbb{P}^4$ be a smooth
quadric hypersurface, and let $h = c_1(\mathcal{O}_{\mathbb{Q}^3}(1))$.
Then $(c_1,c_2,c_3) = (3h, 4 h^2, 2 h^3)$ and $h^3=2$. Letting
$v_1= xh$ and $v_2=y(h^2/2)$ (recall that $h^2/2$ generates
$1$-cycles) for some positive integers $x,y$, and arguing in a
similar way, we get the following result.

\begin{thm} \label{uninflectedoverQ3}
Let $X \subset \mathbb{P}^{N}$ be a fourfold scroll over
$\mathbb{Q}^3$, satisfying our general assumptions. Then
$\Phi_2(X) \not= \emptyset$.
\end{thm}

\medskip
\noindent {\it{Acknowledgements.}} The first author would like to
thank the MiUR of the Italian Government for support received in
the framework of the PRIN ''Algebraic varieties etc.'' (Cofin 2006
and 2008) as well as the University of Milano (FIRST 2006 and
2007) for making this collaboration possible. The second author
wants to thank for the funds supporting this research from the
projects MTM 2007-61124, MTM 2009-06964, and the NILS  Project.

\end{document}